\documentclass[a4paper]{article}

\usepackage[english]{babel}
\usepackage[utf8]{inputenc}
\usepackage[colorinlistoftodos]{todonotes}
\usepackage{amsmath}
\usepackage{amssymb}
\usepackage{amsthm} 
\usepackage{mathtools}
\usepackage{mathrsfs}
\usepackage{siunitx}
\usepackage{algpseudocode}

\usepackage{subfig}

\usepackage{tikz}

\usepackage{tabularx}
\usepackage{booktabs}
\usepackage{multirow}

\usepackage{empheq}

\usepackage{lipsum}

\usepackage[toc,page]{appendix}

\usepackage{rotating}

\theoremstyle{plain}
\theoremstyle{plain}
\theoremstyle{plain}
\theoremstyle{plain}
\theoremstyle{plain}
\theoremstyle{definition}
\theoremstyle{remark}

\everymath{\displaystyle}
\usepackage[a4paper, total={6in, 8in}]{geometry}
\usepackage{csquotes}
\usepackage[backend=bibtex,style=numeric-comp,giveninits=true,doi=false,isbn=false,url=false,eprint=false,maxbibnames=99]{biblatex}

\newif\ifdouble

\newcommand{\papertitle}{Fast and robust parameter estimation with uncertainty quantification for the cardiac function}

\newcommand{\keywordOne}{Cardiac electromechanics}
\newcommand{\keywordTwo}{Machine Learning}
\newcommand{\keywordThree}{Surrogate modeling}
\newcommand{\keywordFour}{Parameter estimation}
\newcommand{\keywordFive}{Uncertainty quantification}

\newcommand{\modEM}{\mathcal{EM}}
\newcommand{\modEMfom}{\mathcal{EM}_{\text{3D}}}
\newcommand{\modEMredANN}{\mathcal{EM}_{\text{ANN}}}
\newcommand{\modEMredemulator}{\mathcal{EM}_{\text{EMULATOR}}}
\newcommand{\modEMzeroD}{\mathcal{EM}_{\text{0D}}}
\newcommand{\modCirc}{\mathcal{C}}

\newcommand{\EMState}{\mathbf{y}}

\newcommand{\CircState}{\mathbf{c}}

\newcommand{\ANNState}{\mathbf{z}}
\newcommand{\ANNStateTilde}{\widetilde{\mathbf{z}}}
\newcommand{\EMRhs}  {\mathcal{L}}
\newcommand{\CircRhs}{\mathbf{f}}
\newcommand{\ANNRhs}{\mathcal{NN}}

\newcommand{\NumANNState} {N_{z}}
\newcommand{\NumANNWeights} {N_{w}}

\newcommand{\NoiseCov}{\boldsymbol{\Sigma}}

\newcommand{\Noisemeas}{\sigma_{\text{meas}}}
\newcommand{\identityMatrix}{\mathbb{I}}

\newcommand{\VLVemFOM}{V_{\mathrm{LV}}^{\mathrm{3D}}}
\newcommand{\VLVemRED}{V_{\mathrm{LV}}^{\mathrm{ANN}}}
\newcommand{\VLVcirc} {V_{\mathrm{LV}}^{\mathrm{0D}}}

\newcommand{\ANNparam}{\mathbf{w}}
\newcommand{\ANNparamTrained}{\widehat{\ANNparam}}

\newcommand{\param}{\boldsymbol{\theta}}
\newcommand{\paramC}{\param_{\mathcal{C}}}
\newcommand{\paramM}{\param_{\mathcal{EM}}}
\newcommand{\paramSpace}{\boldsymbol{\Theta}}

\newcommand{\NumParams} {N_\mathcal{P}}
\newcommand{\NumMomentum} {N_{\rho}}

\newcommand{\paramSpaceC} {\paramSpace_{\mathcal{C}}}
\newcommand{\paramSpaceM} {\paramSpace_{\mathcal{EM}}}
\newcommand{\NumParamC} {N_{\mathcal{C}}}
\newcommand{\NumParamM} {N_{\mathcal{EM}}}

\newcommand{\EMCState}{\mathbf{s}}
\newcommand{\EMCStateInit}{\mathbf{s}_\mathrm{0}}
\newcommand{\adjoint} {\mathbf{a}}
\newcommand{\adjointInit} {\mathbf{a}_\mathrm{T}}

\newcommand{\momentum}{\boldsymbol{\rho}}

\newcommand{\actshape}{\phi}
\newcommand{\tC}{t^\mathrm{contr}}
\newcommand{\tR}{t^\mathrm{rel}}
\newcommand{\TC}{T^\mathrm{contr}}
\newcommand{\TR}{T^\mathrm{rel}}
\newcommand{\tCstar}{t^\mathrm{contr}_{i}}
\newcommand{\tRstar}{t^\mathrm{rel}_{i}}
\newcommand{\TCstar}{T^\mathrm{contr}_{i}}
\newcommand{\TRstar}{T^\mathrm{rel}_{i}}
\newcommand{\THB}{T_\mathrm{HB}}

\newcommand{\Vheart}{V_{\mathrm{heart}}^{\mathrm{tot}}}
\newcommand{\VLA}{V_{\mathrm{LA}}}
\newcommand{\VLV}{V_{\mathrm{LV}}}
\newcommand{\VRA}{V_{\mathrm{RA}}}
\newcommand{\VRV}{V_{\mathrm{RV}}}

\newcommand{\VnLA}{V_{\mathrm{0,LA}}}
\newcommand{\VnLV}{V_{\mathrm{0,LV}}}
\newcommand{\VnRA}{V_{\mathrm{0,RA}}}
\newcommand{\VnRV}{V_{\mathrm{0,RV}}}

\newcommand{\PLA}{p_{\mathrm{LA}}}
\newcommand{\PLV}{p_{\mathrm{LV}}}
\newcommand{\PRA}{p_{\mathrm{RA}}}
\newcommand{\PRV}{p_{\mathrm{RV}}}
\newcommand{\ELA}{E_{\mathrm{LA}}}
\newcommand{\ELV}{E_{\mathrm{LV}}}
\newcommand{\ERA}{E_{\mathrm{RA}}}
\newcommand{\ERV}{E_{\mathrm{RV}}}
\newcommand{\Estar}{E_{\text{i}}}
\newcommand{\EpLA}{E_{\mathrm{LA}}^{\mathrm{pass}}}
\newcommand{\EpLV}{E_{\mathrm{LV}}^{\mathrm{pass}}}
\newcommand{\EpRA}{E_{\mathrm{RA}}^{\mathrm{pass}}}
\newcommand{\EpRV}{E_{\mathrm{RV}}^{\mathrm{pass}}}
\newcommand{\Epstar}{E_{\text{i}}^{\mathrm{pass}}}

\newcommand{\EaRV}{E_{\mathrm{RV}}^{\mathrm{act}}}

\newcommand{\EaMaxLA}{E_{\mathrm{LA}}^{\mathrm{act,max}}}
\newcommand{\EaMaxLV}{E_{\mathrm{LV}}^{\mathrm{act,max}}}
\newcommand{\EaMaxRA}{E_{\mathrm{RA}}^{\mathrm{act,max}}}
\newcommand{\EaMaxRV}{E_{\mathrm{RV}}^{\mathrm{act,max}}}
\newcommand{\EaMaxstar}{E_{\text{i}}^{\mathrm{act,max}}}

\newcommand{\QarSYS}{Q_{\mathrm{AR}}^{\mathrm{SYS}}}
\newcommand{\QarPUL}{Q_{\mathrm{AR}}^{\mathrm{PUL}}}
\newcommand{\QvnSYS}{Q_{\mathrm{VEN}}^{\mathrm{SYS}}}
\newcommand{\QvnPUL}{Q_{\mathrm{VEN}}^{\mathrm{PUL}}}

\newcommand{\CarSYS}{C_{\mathrm{AR}}^{\mathrm{SYS}}}
\newcommand{\CarPUL}{C_{\mathrm{AR}}^{\mathrm{PUL}}}
\newcommand{\CvnSYS}{C_{\mathrm{VEN}}^{\mathrm{SYS}}}
\newcommand{\CvnPUL}{C_{\mathrm{VEN}}^{\mathrm{PUL}}}

\newcommand{\ParSYS}{p_{\mathrm{AR}}^{\mathrm{SYS}}}
\newcommand{\ParPUL}{p_{\mathrm{AR}}^{\mathrm{PUL}}}
\newcommand{\PvnSYS}{p_{\mathrm{VEN}}^{\mathrm{SYS}}}
\newcommand{\PvnPUL}{p_{\mathrm{VEN}}^{\mathrm{PUL}}}

\newcommand{\RarSYS}{R_{\mathrm{AR}}^{\mathrm{SYS}}}
\newcommand{\RarPUL}{R_{\mathrm{AR}}^{\mathrm{PUL}}}
\newcommand{\RvnSYS}{R_{\mathrm{VEN}}^{\mathrm{SYS}}}
\newcommand{\RvnPUL}{R_{\mathrm{VEN}}^{\mathrm{PUL}}}

\newcommand{\LarSYS}{L_{\mathrm{AR}}^{\mathrm{SYS}}}
\newcommand{\LarPUL}{L_{\mathrm{AR}}^{\mathrm{PUL}}}
\newcommand{\LvnSYS}{L_{\mathrm{VEN}}^{\mathrm{SYS}}}
\newcommand{\LvnPUL}{L_{\mathrm{VEN}}^{\mathrm{PUL}}}

\newcommand{\QAV}{Q_{\mathrm{AV}}}
\newcommand{\QMV}{Q_{\mathrm{MV}}}
\newcommand{\QTV}{Q_{\mathrm{TV}}}
\newcommand{\QPV}{Q_{\mathrm{PV}}}
\newcommand{\RAV}{R_{\mathrm{AV}}}
\newcommand{\RMV}{R_{\mathrm{MV}}}
\newcommand{\RTV}{R_{\mathrm{TV}}}
\newcommand{\RPV}{R_{\mathrm{PV}}}
\newcommand{\Rmin}{R_{\mathrm{min}}}
\newcommand{\Rmax}{R_{\mathrm{max}}}

\newcommand{\aXB}{{a_{\text{XB}}}}

\newcommand{\Tall}{\mathcal{T}_{\text{all}}}

\newcommand{\TLV}{\mathcal{T}_{\text{LV}}}

\newcommand{\TLVperturbed}{\mathcal{T}_{\text{LV}}^{\text{perturbed}}}
\graphicspath{{pictures/}}
\doublefalse

\title{{\papertitle}}
\author{Matteo Salvador$^{1,*}$,
        Francesco Regazzoni$^1$,
		Luca Dede'$^1$,
		Alfio Quarteroni$^{1, 2}$}

\date{\footnotesize
	$^1$ MOX-Dipartimento di Matematica, Politecnico di Milano, Milan, Italy \\
	$^2$ \'Ecole Polytechnique F\'ed\'erale de Lausanne, Lausanne, Switzerland (\textit{Professor Emeritus})\\[2ex]
	$^*$ \textit{Corresponding author} (\texttt{matteo1.salvador@polimi.it}) \\
    }

\begin{document}	
	\maketitle

	\begin{abstract}
		Parameter estimation and uncertainty quantification are crucial in computational cardiology, as they enable the construction of digital twins that faithfully replicate the behavior of physical patients.
Robust and efficient mathematical methods must be designed to fit many model parameters starting from a few, possibly non-invasive, noisy observations.
Moreover, the effective clinical translation requires short execution times and a small amount of computational resources.
In the framework of Bayesian statistics, we combine Maximum a Posteriori estimation and Hamiltonian Monte Carlo to find an approximation of model parameters and their posterior distributions.
To reduce the computational effort, we employ an accurate Artificial Neural Network surrogate of 3D cardiac electromechanics model coupled with a 0D cardiocirculatory model.
Fast simulations and minimal memory requirements are achieved by using matrix--free methods, automatic differentiation and automatic vectorization.
Furthermore, we account for the surrogate modeling error and measurement error.
We perform three different \textit{in silico} test cases, ranging from the ventricular function to the entire cardiovascular system, involving whole-heart mechanics, arterial and venous circulation.
The proposed method is robust when high levels of signal-to-noise ratio are present in the quantities of interest in combination with a random initialization of the model parameters in suitable intervals.
As a matter of fact, by employing a single central processing unit on a standard laptop and a few hours of computations, we attain small relative errors for all model parameters and we estimate posterior distributions that contain the true values inside the $90\%$ credibility regions.
With these benefits, our approach meets the requirements for clinical exploitation, while being compliant with Green Computing practices.
	\end{abstract}


	\section{Introduction}
\label{sec:introduction}
Personalization of computational heart models is necessary to better addressing patient--specific pathophysiology and for assisting the clinicians in the decision--making process for medical treatment \cite{Mirams2020, Trayanova2011}.
In this field, sophisticated cell-to-organ level mathematical models comprising systems of nonlinear differential equations, along with efficient and accurate numerical methods, have been developed to properly describe the physical phenomena underlying the cardiac function \cite{Augustin2021, Fedele2022, Gerach2021, Peirlinck2021, Piersanti2022, Regazzoni2022, Strocchi2020Cohort, Strocchi2020Simulating}.

Numerical simulations using these biophysically detailed and anatomically accurate mathematical models call for high computational costs and for a significant amount of computational resources \cite{Regazzoni2022EMROM}.
Imaging techniques are combined with numerical simulations to perform robust parameter estimation in patient--specific cases \cite{Marchesseau2013, Marx2020, Salvador2021, Sermesant2012}.
On the other hand, multi-fidelity models of cardiac electromechanics, deep learning-based models of cardiac mechanics or simplified lumped circulation models are also employed for the same purpose \cite{Cicci2022, Jung2022, Regazzoni2021cardioemulator, Schiavazzi2017}.
All these mathematical tools mainly focus on the ventricular activity of the human heart. 

In this paper, we present a numerical strategy to perform parameter calibration with uncertainty quantification (UQ) by means of a reduced--order model (ROM) of 3D cardiac electromechanics coupled with closed--loop blood circulation \cite{Regazzoni2022EMROM}.
The ROM, which is based on Artificial Neural Networks (ANNs), encodes the dynamics of the pressure-volume relationship obtained from an accurate full--order model (FOM) of the cardiac function \cite{Fedele2022, Piersanti2022, Regazzoni2022}.
Moreover, it allows for real-time numerical simulations on a personal computer while embedding electromechanical parameters of the 3D mathematical model \cite{Regazzoni2022EMROM}.

Parameter estimation is carried out by solving a constrained optimization problem with an efficient adjoint-based method that exploits matrix--free methods, automatic differentiation and automatic vectorization \cite{Chen2019, Manzoni2021}.
Then, we account for the uncertainty coming from possible model and measurement errors. 
Specifically, we employ the Hamiltonian Monte Carlo (HMC) algorithm to perform inverse UQ \cite{Betancourt2013, Homan2014}.

We verify our approach against \textit{in silico} data with different levels of signal--to--noise (SNR) ratio.
We consider several non-invasive time-dependent quantities of interest (QoIs), such as arterial systemic pressure, atrial and ventricular volumes, in order to estimate many model parameters, ranging from cardiac mechanics to cardiovascular hemodynamics.
This can be done in a few hours of total execution time while simply employing one Central Processing Unit (CPU) of a standard laptop.
Our method can therefore be applied to clinical data, where accuracy, robustness and timeliness are certainly essential.

	\section{Mathematical models}
\label{sec:models}

\begin{figure}
    \centering
    \includegraphics[width=0.8\textwidth]{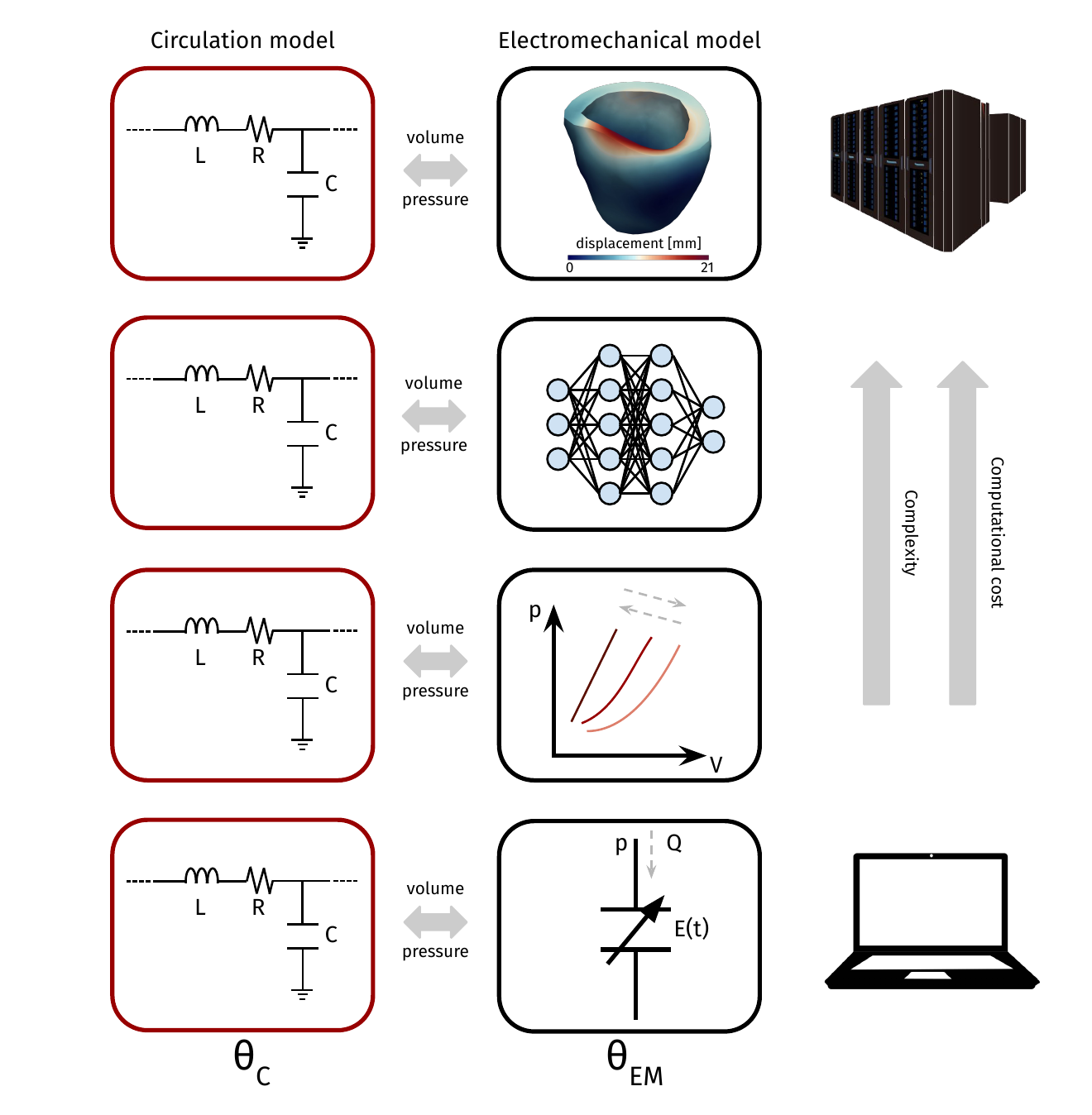}
    \caption{Representation of different $\modEM-\modCirc$ models. A generic circulation model $\modCirc$ may be coupled with biophysically detailed and anatomically accurate $\modEMfom$, ANN-based $\modEMredANN$, emulator-based $\modEMredemulator$ or 0D $\modEMzeroD$ electromechanical models according to the requirements in terms of accuracy, computational efficiency and memory storage.}
    \label{fig:EMcirculation}
\end{figure}

We display in Figure~\ref{fig:EMcirculation} several mathematical models for the cardiac function featuring a different degree of physical accuracy and computational complexity.
These mathematical models can be regarded as a FOM for cardiac electromechanics $\modEMfom$ and three different ROMs.
Although model $\modEMfom$ presents a high level of biophysical accuracy, it is also associated to high performance computing and significant computational costs.
This motivates the use of ROMs, which are computationally cheaper than the corresponding FOM in terms of execution time and computer resources.
Moreover, they do not significantly compromise the accuracy of the FOM.
These ROMs simulate the pressure-volume relationship of one or multiple cardiac chambers and can be all employed for fast and robust parameter estimation. We will use the following taxonomy:
\begin{itemize}
\item $\modEMredANN$: ANN based surrogate models of cardiac electromechanics, built as black-box from a collection of pre-computed numerical simulations through a data-driven approach \cite{Regazzoni2022EMROM};
\item $\modEMredemulator$: parametric emulators of cardiac electromechanics built with a grey-box approach, that is by fitting a priori defined physics-inspired curves from data obtained by means of numerical simulations \cite{Regazzoni2021cardioemulator};
\item $\modEMzeroD$: fully 0D electromechanical models, i.e. time-varying elastance models assuming a linear relationship between pressure and volume \cite{Hirschvogel2017, Regazzoni2022}.
\end{itemize}
The aforementioned mathematical models are coupled with a generic circulation model $\modCirc$ of the remaining part of the cardiovascular system by exchanging pressures and volumes.
For the sake of simplicity, in this paper we consider a ROM built with the $\modEMredANN$ model for the left ventricle (LV) only.

\subsection{3D electromechanical model}
\label{sec:models:electromechanics}

We model the electromechanical activity of the LV by means of a set of differential equations.
In compact form we can write the model as \cite{Regazzoni2021cardioemulator, Regazzoni2022EMROM}:
\begin{equation} \label{eqn:model_3D}
    \left\{
    \begin{aligned}
        \frac{\partial \EMState(t)}{\partial t} &= \EMRhs  (\EMState(t)  , \PLV(t), t; \paramM) && \text{for } t \in (0,T],\\
        \EMState  (0) &= \EMState_0.   && \\
    \end{aligned}
    \right.
\end{equation}
The nonlinear differential operator $\EMRhs$ encodes the differential equations and boundary conditions, 
while the state vector $\EMState(t)$ contains several variables associated with the cardiac function, such as the action potential, intracellular calcium concentration, sarcomere length and mechanical deformation.
$\EMState_0$ represents the initial condition.
$\PLV(t)$ indicates the endocardial pressure of the LV.
We introduce $\paramSpaceM \subseteq \mathbb{R}^{\NumParamM}$, that is the parameter space, being $\NumParamM$ the number of parameters, and we denote as $\paramM$ the model parameters for cardiac electromechanics such that $\paramM \in \paramSpaceM$. Examples of electromechanical parameters are electric conductivities, fibers direction, contractility and passive stiffness of the myocardium.

From now on, we label Equation~\eqref{eqn:model_3D} as $\modEMfom$, which requires a closure relationship to determine the pressure $\PLV(t)$.
Indeed, we couple $\modEMfom$ with a generic circulation model of the cardiovascular system, henceforth denoted as $\modCirc$ \cite{Augustin2021, Peirlinck2021, Piersanti2022, Regazzoni2022}.
The coupled problem reads \cite{Regazzoni2021cardioemulator, Regazzoni2022EMROM}:
\begin{equation} \label{eqn:model_3D-0D}
    \left\{
    \begin{aligned}
        \frac{\partial \EMState(t)}{\partial t} &= \EMRhs  (\EMState(t)  , \PLV(t), t; \paramM) && \text{for } t \in (0,T],\\
        \frac{d \CircState(t)}{dt} &= \CircRhs(\CircState(t), \PLV(t), t; \paramC) && \text{for } t \in (0,T],\\
        \VLVemFOM(\EMState(t)) &= \VLVcirc(\CircState(t)) && \text{for } t \in (0,T],\\
        \EMState  (0) &= \EMState_0,   && \\
        \CircState(0) &= \CircState_0. && \\
    \end{aligned}
    \right.
\end{equation}
The state variables $\CircState(t)$ of the cardiocirculatory model contain pressures, volumes and fluxes in different compartments of the vascular network, while $\paramC \in \paramSpaceC \subseteq \mathbb{R}^{\NumParamC}$ is a vector of $\NumParamC$ parameters that contains, for instance, resistances, conductances or elastances.
The volumetric constraint $\VLVemFOM(\EMState(t)) = \VLVcirc(\CircState(t))$ between $\modEMfom$ and $\modCirc$ models is enforced by means of $\PLV(t)$ \cite{Regazzoni2022}, which acts numerically as a Lagrange multiplier \cite{Regazzoni2022}.
An alternative approach would be to adapt the closure relationships to the different phases of the heartbeat \cite{levrero2020sensitivity, Salvador2020}, e.g. by using preload and afterload windkessel models \cite{Westerhof2009}. 

The surrogate models that will be introduced in the next sections are built from the 3D-0D closed-loop mathematical model $\modEMfom-\modCirc$ presented in \cite{Piersanti2022, Regazzoni2022}.
We consider a LV obtained from the Zygote 3D human heart model \cite{zygote}, endowed with a fiber architecture generated by means of suitable Laplace-Dirichlet-Rule-Based methods \cite{Piersanti2020}.
Nevertheless, our mathematical and numerical models directly generalize to patient-specific geometries \cite{Salvador2021}.
For cardiac electrophysiology, we employ the monodomain equation \cite{collifranzone2014book} coupled with the ten Tusscher-Panfilov ionic model \cite{ten2006alternans}.
We model mechanical activation in the active stress formulation by means of the biophysically detailed and anatomically accurate RDQ20-MF model \cite{regazzoni2020biophysically}.
We consider the Guccione constitutive law in a quasi-incompressible regime for passive mechanics \cite{usyk2002computational}.
We adopt spring-damper Robin boundary conditions at the epicardium to account for the presence of the pericardial sac \cite{Gerbi2018,Pfaller2019}.
We prescribe energy-consistent boundary conditions to model the interaction with the part of the myocardium beyond the artificial ventricular base \cite{regazzoni2019mor-sarcomeres}.
Blood circulation over the whole cardiovascular system is modeled by using a closed-loop mathematical model $\modCirc$ proposed in \cite{Hirschvogel2017, Regazzoni2022}. For the sake of completeness, we report the equations of the 0D cardiocirculatory model in Appendix~\ref{app:circulation}.

\subsection{Artificial neural network based reduced-order model}
\label{sec:models:electromechanics:ANN}

Following the ANN-based ROM introduced in \cite{Regazzoni2022EMROM}, we build a set of ordinary differential equations (ODEs), whose right hand side is represented by an ANN, that learns the pressure-volume dynamics of the 3D cardiac electromechanical model $\modEMfom$ reported in Equation~\eqref{eqn:model_3D}.
In this framework, the ANN-based ROM $\modEMredANN$ for the LV reads: 
\begin{equation} \label{eqn:model_ANN}
    \left\{
    \begin{aligned}
        \frac{d \ANNState(t)}{d t} &= \ANNRhs\left(
            \ANNState(t),
            \PLV(t),
            \cos(\tfrac{2 \pi t}{\THB}),
            \sin(\tfrac{2 \pi t}{\THB}),
            \paramM;
            \ANNparamTrained
        \right)
        && \text{for } t \in (0,T],\\
        \ANNState  (0) &= \ANNState_0.  && \\
    \end{aligned}
    \right.
\end{equation}
The fully connected feedforward ANN is defined by $\ANNRhs \colon \mathbb{R}^{\NumANNState + \NumParamM + 3} \to \mathbb{R}^{\NumANNState}$ and $\ANNState(t) \in \mathbb{R}^{\NumANNState}$ represents the reduced state vector.
The ANN receives $\NumANNState$ state variables, $\NumParamM$ scalar parameters, pressure $\PLV$, and two periodic conditions in time as input.
Indeed, the $\cos({2 \pi t}/{\THB})$ and $\sin({2 \pi t}/{\THB})$ terms account for the periodicity $\THB$ of the heartbeat, thus allowing for arbitrarily long real-time numerical simulations of cardiac electromechanics \cite{Regazzoni2022EMROM}.
Finally, weights and biases of the ANN are encoded in $\ANNparamTrained \in \mathbb{R}^{\NumANNWeights}$.

We train the ANN to enable the predicted LV volume to coincide with the first state variable, i.e. $\ANNState(t) = [\VLVemRED(\ANNState(t)), \ANNStateTilde(t)]^T$ with $\ANNStateTilde(t) \in \mathbb{R}^{\NumANNState-1}$.
In this manner, by following the notation introduced in \cite{regazzoni2019modellearning}, we consider an output-inside-the-state approach.
The remaining components of vector $\ANNState(t)$, that is $\ANNStateTilde(t)$, are latent variables without immediate physical interpretation.
For all the details regarding the model training strategy of this ANN-based ROM we refer to \cite{Regazzoni2022EMROM}.

The coupled $\modEMredANN-\modCirc$ model gives rise to a differential-algebraic system of equations (DAEs) that reads:
\begin{equation} \label{eqn:model_ANN_circulation}
    \left\{
    \begin{aligned}
        \frac{d \ANNState(t)}{d t} &= \ANNRhs\left(
            \ANNState(t),
            \PLV(t),
            \cos(\tfrac{2 \pi t}{\THB}),
            \sin(\tfrac{2 \pi t}{\THB}),
            \paramM;
            \ANNparamTrained
        \right)
        && \text{for } t \in (0,T], \\
        \frac{d \CircState(t)}{dt} &= \CircRhs(\CircState(t), \PLV(t), t; \paramC) && \text{for } t \in (0,T], \\
        \VLVcirc(\CircState(t)) &= \VLVemRED(\ANNState(t)) && \text{for } t \in (0,T], \\
        \ANNState (0) &= \ANNState_0,   && \\
        \CircState(0) &= \CircState_0. && \\
    \end{aligned}
    \right.
\end{equation}
In order to perform parameter estimation by employing an adjoint sensitivity method \cite{Chen2019}, we need to convert this DAEs into a system of ODEs:
\begin{equation} \label{eqn:model_ANN_circulation}
    \left\{
    \begin{aligned}
        \frac{d \ANNState(t)}{d t} &= \ANNRhs\left(
            \ANNState(t),
            \PLV(t),
            \cos(\tfrac{2 \pi t}{\THB}),
            \sin(\tfrac{2 \pi t}{\THB}),
            \paramM;
            \ANNparamTrained
        \right)
        && \text{for } t \in (0,T], \\
        \frac{d \CircState(t)}{dt} &= \CircRhs(\CircState(t), \PLV(t), t; \paramC) && \text{for } t \in (0,T], \\
        \varepsilon \frac{d \PLV(t)}{d t} &= \VLVemRED(\ANNState(t)) - \VLVcirc(\CircState(t)) && \text{for } t \in (0,T], \\
        \ANNState (0) &= \ANNState_0,   && \\
        \CircState(0) &= \CircState_0. && \\
    \end{aligned}
    \right.
\end{equation}
$\varepsilon$ is a small penalization parameter; for $\varepsilon \to 0$, we recover the original constraint $\VLVcirc(\CircState(t)) = \VLVemRED(\ANNState(t))$.

	\section{Mathematical methods}
\label{sec:methods}

We define an optimal control problem and we describe the mathematical techniques that we use to perform parameter estimation with inverse UQ.

\subsection{Optimization problem}
\label{sec:methods:OCP}
Let $\param = [\paramM, \paramC]^T \in \paramSpace \subset \mathbb{R}^{\NumParams}$ be a set of parameters for the $\modEMredANN-\modCirc$ model.
Then, we want to solve a constrained optimization problem \cite{Manzoni2021}:
\begin{equation}
\min_{\param \in \paramSpace} J(\param)
\label{eqn: OCP}
\end{equation}
with the following cost functional:
\begin{equation}
\begin{split}
J(\param) &= \sum_{\text{i} \in \{ \text{LA, LV, RA, RV, AR-SYS} \}} \alpha_\text{i} \dfrac{||p_\text{i}(t;\param) - \hat{p}_\text{i}(t)||_{\text{L}^2}^2}{\mu_{p_\text{i}^2}} \\
& + \sum_{\text{i} \in \{ \text{LA, LV, RA, RV, AR-SYS} \}} \beta_\text{i} \dfrac{|p_\text{i}^\text{max}(\param) - \hat{p}_\text{i}^\text{max}|}{\mu_{p_\text{i}^2}}
^2 \\
& + \sum_{\text{i} \in \{ \text{LA, LV, RA, RV, AR-SYS} \}} \gamma_\text{i} \dfrac{|p_\text{i}^\text{min}(\param) - \hat{p}_\text{i}^\text{min}|^2}{\mu_{p_\text{i}^2}} \\
& + \sum_{\text{i} \in \{ \text{LA, LV, RA, RV} \}} \delta_\text{i} \dfrac{||V_\text{i}(t;\param) - \hat{V}_\text{i}(t)||_{\text{L}^2}^2}{\mu_{V_\text{i}^2}} \\
& + \sum_{\text{i} \in \{ \text{LA, LV, RA, RV} \}} \epsilon_\text{i} \dfrac{|V_\text{i}^\text{max}(\param) - \hat{V}_\text{i}^\text{max}|^2}{\mu_{V_\text{i}^2}} \\
& + \sum_{\text{i} \in \{ \text{LA, LV, RA, RV} \}} \zeta_\text{i} \dfrac{|V_\text{i}^\text{min}(\param) - \hat{V}_\text{i}^\text{min}|^2}{\mu_{V_\text{i}^2}} \\
& + \sum_{\text{i} \in \{ \text{LA, LV, RA, RV} \}} \eta_\text{i} \dfrac{|SV_\text{i}(\param) - \hat{SV}_\text{i}|^2}{\mu_{V_\text{i}^2}},
\end{split}
\label{eqn: cost_functional}
\end{equation}
where $t \in [T - \THB, T]$, for $T$ sufficiently large to reach a limit cycle.
Problem~\eqref{eqn: OCP} is constrained in the sense that all the pressures $p_\text{i}$ and volumes $V_\text{i}$ appearing in~\eqref{eqn: cost_functional} are solutions of the $\modEMredANN-\modCirc$ model.
On the other hand, the corresponding observations $\hat{p}_\text{i}$ and $\hat{V}_\text{i}$ are general, i.e. they may come from in silico numerical simulations or clinical data. 
Coefficients $\alpha_\text{i}$ to $\eta_\text{i}$ weigh the different contributions, which involve pressure and volume traces over time as well as pointwise quantities, such as maximum, minimum and stroke volume (SV) values. In particular, these quantities are defined as follows:
\begin{equation*}
    \begin{array}{l}
        p_\text{i}^\text{max}(\param) = \max_{t \in [T - \THB, T]} p_\text{i}(t;\param) \;\;\; \text{for} \;\;\; \text{i} \in \{ \text{LA, LV, RA, RV, AR-SYS} \} \\[10pt]
        V_\text{i}^\text{max}(\param) = \max_{t \in [T - \THB, T]} V_\text{i}(t;\param) \;\;\; \text{for} \;\;\; \text{i} \in \{ \text{LA, LV, RA, RV} \} \\[10pt]
        p_\text{i}^\text{min}(\param) = \min_{t \in [T - \THB, T]} p_\text{i}(t;\param) \;\;\; \text{for} \;\;\; \text{i} \in \{ \text{LA, LV, RA, RV, AR-SYS} \} \\[10pt]
        V_\text{i}^\text{min}(\param) = \min_{t \in [T - \THB, T]} V_\text{i}(t;\param) \;\;\; \text{for} \;\;\; \text{i} \in \{ \text{LA, LV, RA, RV} \} \\[10pt]
        SV_\text{i}(\param) = V_\text{i}^\text{max}(\param) - V_\text{i}^\text{min}(\param) \;\;\; \text{for} \;\;\; \text{i} \in \{ \text{LA, LV, RA, RV} \}
    \end{array}
\end{equation*}
We also introduce two normalization terms, namely $\mu_{p_\text{i}^2}$ and $\mu_{V_\text{i}^2}$, which are obtained by averaging the squared values of pressure and volume traces over time, for $t \in [T - \THB, T]$.

For the sake of simplicity, in this paper we limit ourselves to QoIs regarding the four cardiac chambers and arterial systemic network.
Nevertheless, fluxes accross the four cardiac valves or pressures and fluxes related to the pulmonary circulation can be seamlessly added in the cost functional.

\subsection{Parameter estimation}
\label{sec:methods:PE}

\begin{figure}
    \centering
    \includegraphics[width=1.0\textwidth]{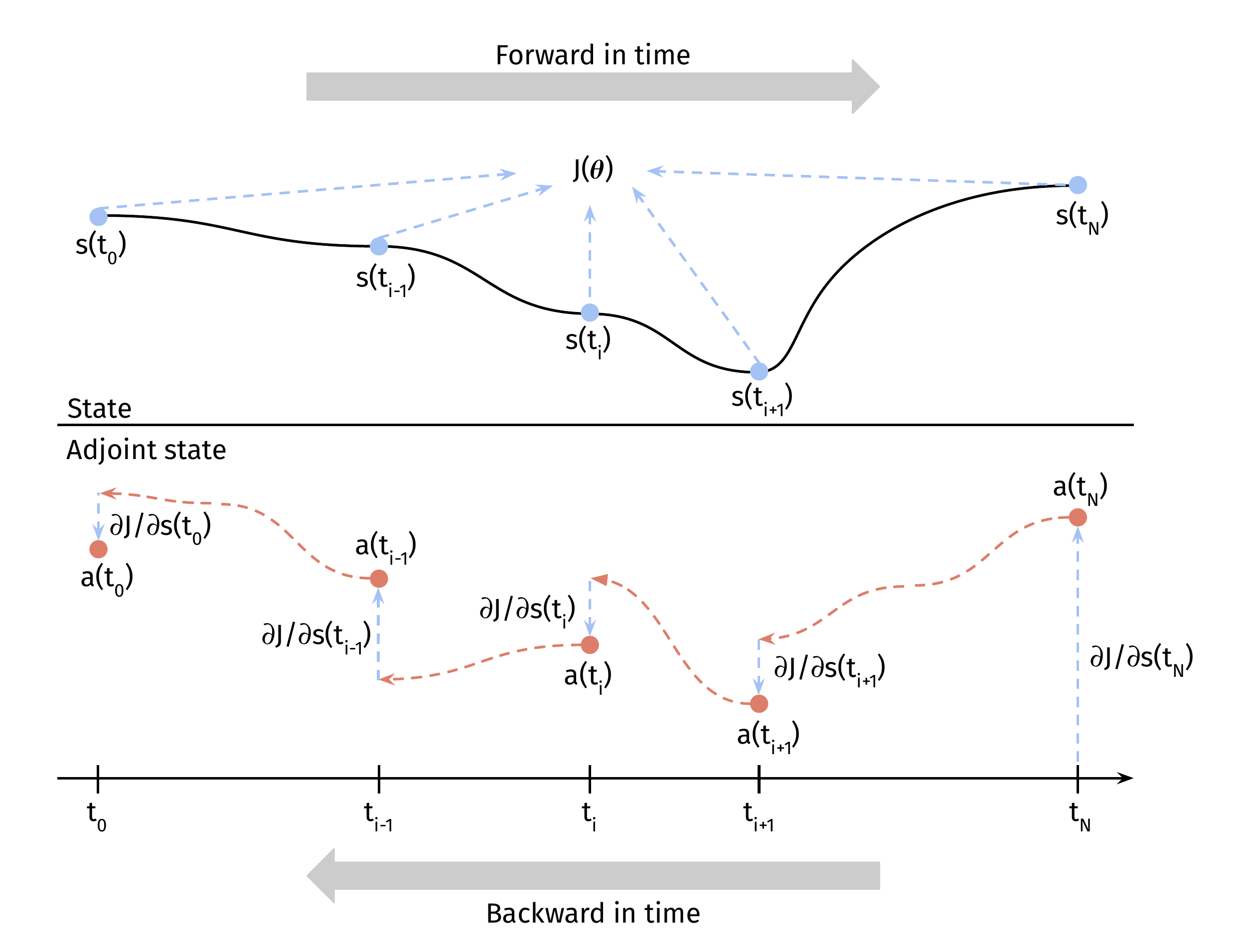}
    \caption{Sketch of the adjoint-based optimization algorithm applied on a single ODE. At each iteration, the ODE is solved forward in time and the cost functional $J(\param)$ is computed. Then, reverse-mode differentiation is used to solve the adjoint-based ODE backward in time and to compute the sensitivities of the cost functional with respect to the state variable. The time step is chosen according to an adaptive procedure.}
    \label{fig:AdjointBasedOpt}
\end{figure}

Equation~\eqref{eqn:model_ANN_circulation} can be compactly rewritten as follows:
\begin{equation} 
\begin{cases}
\dfrac{d \EMCState(t)}{d t} = \mathbf{g}(t, \EMCState(t); \param) & \qquad \text{for } t \in (0, T], \\
\EMCState(0) = \EMCStateInit,
\end{cases}
\label{eqn: forward_problem}
\end{equation}
where $\EMCState(t)$ may contain (reduced) state variables related to cardiac electromechanics and circulation variables, while $\mathbf{g}(t, \EMCState(t); \param)$ represents the $\modEMredANN-\modCirc$ coupled problem.
In the optimization process, at each iteration, we have to compute $\tfrac{d J}{d \param}$ to minimize the loss function $J = J(\param)$.
The first step consists in solving Equation~\eqref{eqn: forward_problem} forward in time, i.e. for $t \in (0, T]$, using an ODE numerical solver.
Then, we introduce the adjoint state variables $\adjoint(t)$ and we solve the following adjoint ODE system backwards in time by means of reverse-mode differentiation \cite{Chen2019}:
\begin{equation} 
\begin{cases}
\dfrac{d \adjoint(t)}{d t} = -\adjoint(t)^T \dfrac{\partial \mathbf{g}}{\partial \EMCState}(t, \EMCState(t); \param) & \qquad \text{for } t \in (T, 0], \\
\adjoint(T) = \adjointInit.
\end{cases}
\label{eqn: adjoint_problem}
\end{equation}
Finally, the gradient of $J$ with respect to $\param$ reads \cite{Chen2019}:
\begin{equation}
\frac{d J}{d \param} = - \int_T^0 \adjoint(t)^T \frac{\partial \mathbf{g}}{\partial \param}(t, \EMCState(t); \param) dt
\label{eqn: loss_gradient}
\end{equation}
The vector-jacobian products $\adjoint(t)^T \tfrac{\partial \mathbf{g}}{\partial \EMCState}(t, \EMCState(t); \param)$ and $\adjoint(t)^T \tfrac{\partial \mathbf{g}}{\partial \param}(t, \EMCState(t); \param)$ are efficiently computed in a matrix-free fashion, accounting for small memory requirements, by also exploiting automatic differentiation and automatic vectorization.
For the sake of clarity, the whole optimization process is depicted in Figure~\ref{fig:AdjointBasedOpt}.

The optimal control problem~\eqref{eqn: OCP} is solved by employing the Limited-memory Broyden-Fletcher-Goldfarb-Shanno (L-BFGS) algorithm \cite{Liu1989}.
We use the Dormand-Prince ODE integrator with an adaptive time step to solve the forward and backward problems \cite{Dormand1980}.
Further details about the software library are provided in Section~\ref{sec:methods:software}.

We remark that, as we are dealing with constrained optimization and bounded intervals for $\param \in \paramSpace$, minimizing $J(\param)$ is equivalent, from a probabilistic standpoint, to find the Maximum a Posteriori (MAP) estimation.
This means that, given a uniform prior distribution on the parameters $\mathbb{P}(\theta_\text{i}) \sim\ U(a_\text{i}, b_\text{i})$ for $i = 1, \dots, \NumParams$, the MAP estimation $\param_\text{MAP}$ reads \cite{Murphy2012}:
\begin{equation}
\param_\text{MAP} = \underset{\param}{\operatorname{arg max}} \; \mathbb{P}(\param | \mathbf{x}),
\label{eqn: MAP}
\end{equation}
where $\mathbb{P}(\param | \mathbf{x})$ denotes the posterior distribution over the observations $\mathbf{x}$.
We remark that this equivalence holds under the hypothesis of Gaussian measurement error.

\subsection{Uncertainty quantification}
\label{sec:methods:UQ}

We perform inverse UQ by employing HMC \cite{Betancourt2013}, namely a Markov Chain Monte Carlo (MCMC) method that aims at finding an approximation of the posterior distribution $\mathbb{P}(\param | \mathbf{x})$. We now briefly recall the principles of HMC.

We introduce a vector of auxiliary momentum variables $\momentum \in \mathbb{R}^{\NumMomentum}$, along with the following joint distribution \cite{Betancourt2013}:
\begin{equation*}
\mathbb{P}(\momentum, \param) = \mathbb{P}(\momentum | \param) \mathbb{P}(\param),
\end{equation*}
where $\mathbb{P}(\momentum | \param)$ represents the conditional probability distribution of $\momentum$ given $\param$, while $\mathbb{P}(\param)$ defines the prior probability distribution with respect to $\param$.
Then, we define the Hamiltonian function as follows \cite{Betancourt2013}:
\begin{equation*}
\mathbb{H}(\momentum, \param) = -\log \mathbb{P}(\momentum, \param) = -\log \mathbb{P}(\momentum | \param) -\log \mathbb{P}(\param) = \mathbb{K}(\momentum | \param) + \mathbb{U}(\param),
\end{equation*}
where $\mathbb{K}(\momentum | \param)=-\log \mathbb{P}(\momentum | \param)$ plays the role of a kinetic energy, while $\mathbb{U}(\param)=-\log \mathbb{P}(\param)$ is a potential energy.
Finally, we solve a coupled system of differential equations in $(\param, \momentum)$ to possibly advance the value of the parameters vector $\param=\param(t)$ from its current state:
\begin{equation} 
\begin{cases}
\dfrac{d \param}{d t} = \dfrac{\partial \mathbb{H}}{\partial \momentum} & \qquad \text{for } t \in (0, \overline{T}], \\
\dfrac{d \momentum}{d t} = -\dfrac{\partial \mathbb{H}}{\partial \param} & \qquad \text{for } t \in (0, \overline{T}],
\end{cases}
\label{eqn: HMC}
\end{equation}
by means of the leapfrog integrator \cite{Betancourt2013}, which requires a suitable choice of the time step $\overline{\Delta t}$ and final time $\overline{T}$.
We remark that here $t$ is a fictitious time related to the evolution of the trajectories that explore new possible values for $\param$.
We employ the No-U-Turn Sampler (NUTS) extension of HMC so that the number of steps is automatically adapted by the algorithm \cite{Homan2014}.
To summarize, UQ is performed by repeatedly following these steps \cite{Betancourt2013}:
\begin{enumerate}
\item a new sample $\momentum^0$ is drawn from a zero-mean Gaussian distribution for the momentum variables, that is $\momentum \sim\ \mathcal{N}(\mathbf{0}, \mathbf{\Sigma})$;
\item starting from $[\param^0, \momentum^0]^T$, solve Equation~\eqref{eqn: HMC} with the leapfrog numerical scheme and NUTS;
\item negate the momentum variables and define a proposed state $[\param^*, \momentum^*]^T$;
\item the proposed state $[\param^*, \momentum^*]^T$ is accepted as the next state using a Metropolis-Hastings update with probability $\min (1, \text{exp}(\mathbb{H}(\momentum^0, \param^0) - \mathbb{H}(\momentum^*, \param^*)))$. On the other hand, if the proposal is not accepted, $\param^0$ is used again to initialize the next state of HMC.
\end{enumerate}
Once a certain (a priori fixed) number of iterations are made, the accepted proposals for $\param$, along with the corresponding probability values, are employed to define an approximated posterior distribution for $\mathbb{P}(\param | \mathbf{x})$.

We fix $\overline{\Delta t}=10^{-3}$.
We perform 750 iterations for all the test cases.
Among them, the first 250 iterations consist of an initial warmup period and are not retained to estimate the posterior distribution.
We initialize the NUTS sampler by considering $\mathbb{P}(\param) \sim\ U(\param_\text{MAP} - \iota \param_\text{MAP}, \param_\text{MAP} + \iota \param_\text{MAP})$, being $\iota=0.1$ a suitable parameter to define a uniform prior distribution around $\param_\text{MAP}$.

We account for both the measurement error and the ROM approximation error in the QoIs of Equation~\eqref{eqn: cost_functional} ($\hat{p}_\text{i}(t)$ and $\hat{V}_\text{i}(t)$).
Indeed, pressures and volumes are intrinsically affected by the sensitivity of the instrument that is used to perform the measurements and by the environment in which these values are extracted.
Furthermore, passing from the FOM to the ANN-based ROM introduces an approximation error that must be incorporated in UQ.
We only focus on non-invasive QoIs, i.e. systemic arterial pressure and volumes of the four cardiac chambers.
We consider normal distributions around the estimated values of these QoIs, i.e. $\mathcal{N}(V_\text{i}(t;\param), \Noisemeas^2 \identityMatrix+\NoiseCov_{\text{ROM,i}})$ for $\text{i} \in \{ \text{LA, LV, RA, RV} \}$ and $\mathcal{N}(p_\text{AR-SYS}(t;\param), \Noisemeas^2 \identityMatrix+\NoiseCov_{\text{ROM,AR-SYS}})$.
To model the ROM approximation error, given the time-dependant (correlated) QoIs that we consider in this paper, we introduce a zero-mean Gaussian process $\mathcal{G}\mathcal{P}(\boldsymbol{0}, k(\boldsymbol{t}, \boldsymbol{t}'))$, where $k(\boldsymbol{t}, \boldsymbol{t}')=\sigma^2 \exp \left( \tfrac{-||\boldsymbol{t}-\boldsymbol{t}'||^2}{2 \lambda^2} \right)$ is the exponentiated quadratic kernel \cite{Rasmussen2005}.
Amplitude $\sigma$ is indipendently computed for all the relevant pressure and volume time traces as the root-mean-square errors of the ROM with respect to the FOM \cite{Regazzoni2022EMROM}.
This leads to $\sigma_{\text{LA}}=0.79 \; \si{\milli\liter}$, $\sigma_{\text{LV}}=0.80 \; \si{\milli\liter}$, $\sigma_{\text{RA}}=0.07 \; \si{\milli\liter}$, $\sigma_{\text{RV}}=0.12 \; \si{\milli\liter}$ and $\sigma_{\text{AR-SYS}}=0.54 \; \si{\mmHg}$.
On the other hand, the correlation length $\lambda$ is estimated by minimizing the negative log likelihood of the observed ROM approximation errors.
We employ 1'000 Adam iterations to optimize this kernel hyperparameter \cite{Kingma2014}.
This overall leads to similar correlation lengths for all the interested QoIs. For this reason, we fix a unique value of $\lambda=0.04$.
Once this tuning process ends, the kernel function can be used to generate different full covariance matrices:
\begin{equation*}
\NoiseCov_{\text{ROM,i}}(t_\text{j}, t_\text{k}) = \sigma_{\text{i}}^2 \exp \left[ \dfrac{-(t_\text{j}-t_\text{k})^2}{2 \lambda^2} \right] \;\;\; \text{for} \;\;\; \text{i} \in \{ \text{LA, LV, RA, RV, AR-SYS} \},
\end{equation*}
where $t_\text{j}$ and $t_\text{k}$ represent discrete time values in the time interval $[T - \THB, T]$.
With respect to vanilla Monte Carlo integration, where multiple chains are usually required to achieve proper convergence to the posterior distribution \cite{Regazzoni2022EMROM}, here we will always run a single chain, as this is sufficient to provide meaningful results.
This is also motivated by the suitable initialization of the parameters, which is related to the MAP estimation $\param_\text{MAP}$.
We declare convergence when the Gelman-Rubin diagnostic provides a value less than 1.1 for all the parameters $\param$ and there are no divergent transitions \cite{Brooks1998,Gelman1992}.

\begin{figure}
    \centering
    \includegraphics[width=1.0\textwidth]{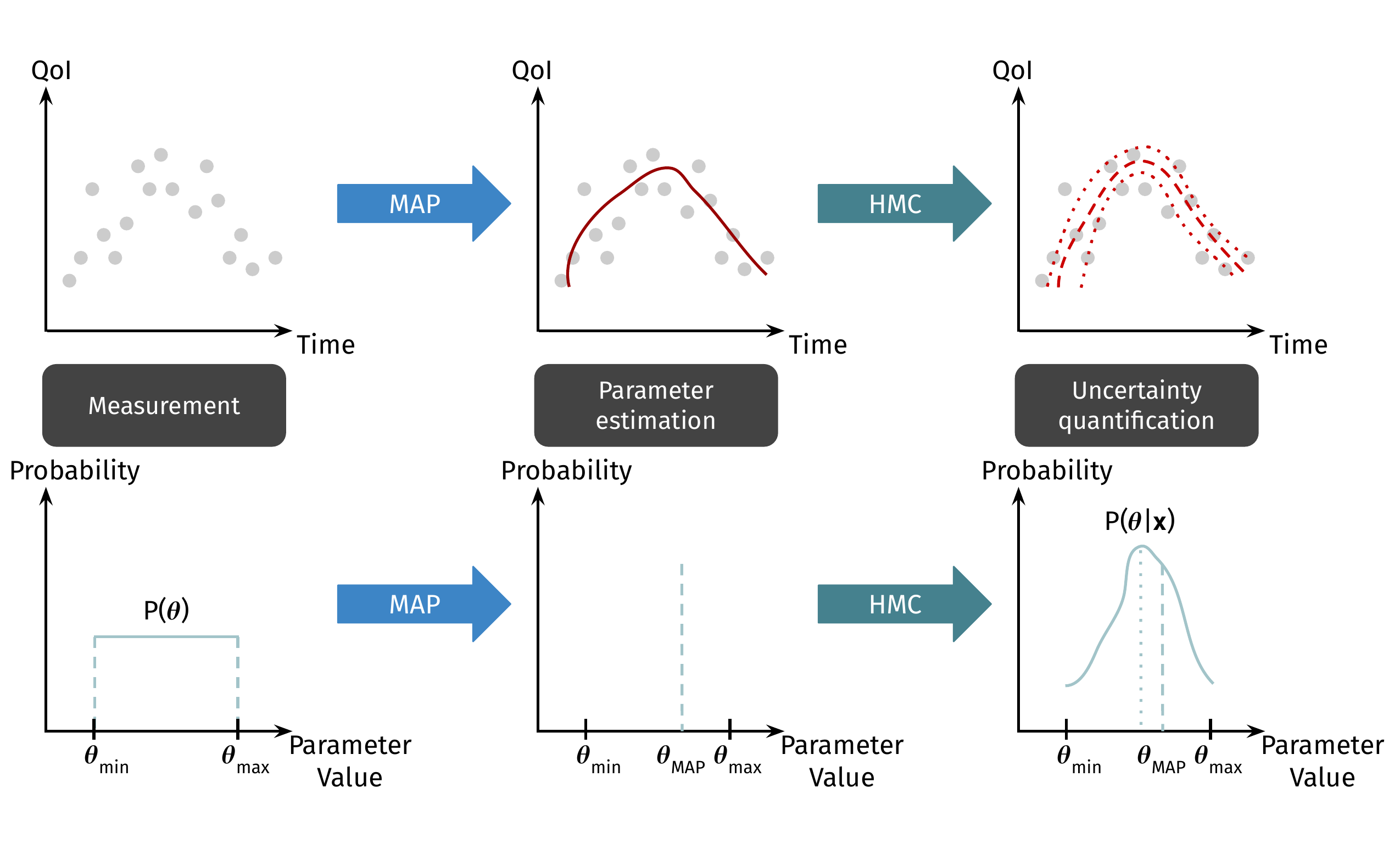}
    \caption{Robust parameter estimation for the cardiac function. Given a generic measurement over time and a uniform prior distribution $\mathbb{P}(\theta)$ of a parameter $\theta$, we compute the MAP estimation $\theta_\text{MAP}$, whose observation fits the (noisy) data. Then, by means of HMC, we compute the posterior distribution $\mathbb{P}(\theta | \mathbf{x})$ and we perform UQ by including both the measurement error and the ROM approximation error.}
    \label{fig:MAPHMC}
\end{figure}

In Figure~\ref{fig:MAPHMC}, we show a synthesis of the whole procedure to perform estimation with uncertainty for a single parameter given a generic QoI of the cardiac function.

\subsection{Software libraries and hardware}
\label{sec:methods:software}
We perform the electromechanical simulations by means of the \texttt{life\textsuperscript{x}} \texttt{C++} library\footnote{\url{https://lifex.gitlab.io}} developed within the iHEART project.
We train the LV ANN-based model with the \texttt{model-learning}\footnote{\url{https://github.com/FrancescoRegazzoni/model-learning}} open source MATLAB library \cite{regazzoni2019modellearning}.
In particular, we exploit the ROM of cardiac electromechanics proposed in \cite{Regazzoni2022EMROM}.
We rely on the \texttt{JAX} open source Python library\footnote{\url{https://github.com/google/jax}} for adjoint based parameter estimation \cite{jax2018github}.
Finally, we carry out HMC based Bayesian inverse UQ by using the \texttt{NumPyro}\footnote{\url{https://github.com/pyro-ppl/numpyro}} open source Python library \cite{phan2019}.
To employ the ANN-based models trained with the \texttt{model-learning} MATLAB library within the Python environments of \texttt{JAX} and \texttt{NumPyro}, we exploit \texttt{pyModelLearning}, a Python wrapper for the \texttt{model-learning} library.

Parameter estimation with UQ has been carried out by using a laptop endowed with 2 cores (one Intel Core i7-7500U CPU, 2.70 GHz) and 16 GB of RAM.
Nevertheless, we run all the test cases by considering 1 core, that is serial execution.

This manuscript is accompanied by \url{https://github.com/MatteoSalvador/cardioEM-MAP}, a public repository containing the codes that allow to reproduce the MAP estimation for all test cases.


	\section{Results}
\label{sec:results}

We design several in silico test cases to show the robustness and flexibility of our approach driven by a combined use of MAP estimation and HMC.
We show that complex scenarios involving many parameters and QoIs can be handled with a small amount of computational resources, i.e. a single core of a standard computer, along with relatively small computational times, that is a few hours.

First, we identify the most important parameters of the electromechanical model by means of global sensitivity analysis (Section~\ref{sec:GSA}).
Then, we formulate the following test cases:

\begin{itemize}
\item $\TLV$: we provide $\ParSYS$ and $\VLV$ time traces and we estimate parameters for the LV and systemic circulation by means of the $\modEMredANN-\modCirc$ model.
\item $\Tall$: we consider the evolution of $\ParSYS$, $\VLA$, $\VLV$, $\VRA$, $\VRV$ over time and we estimate several parameters for the heart and peripheral circulation by means of the $\modEMredANN-\modCirc$ model;
\item $\TLVperturbed$: same as $\TLV$, where we consider perturbed values of relevant parameters in the cardiovascular system with respect to the ground truth. Indeed, in realistic applications, it is uncommon to know a priori all the parameter values that are not related to LV and systemic circulation.
\end{itemize}

We weigh the different contributions of the cost functional (Equation~\ref{eqn: cost_functional}) in the same manner for all the test cases, i.e. all the coefficients ranging from $\alpha_\text{i}$ to $\eta_\text{i}$ have been set to 1 or 0 accordingly.

We mimic the presence of measurement errors for both MAP estimation and HMC by applying different levels of Gaussian noise $\mathcal{N}(\boldsymbol{0}; \Noisemeas^2 \identityMatrix)$ to the QoIs in the different time points.
In particular, we consider $\Noisemeas=0 \; \si{\mmHg}$, $\Noisemeas=1 \; \si{\mmHg}$, $\Noisemeas=10 \; \si{\mmHg}$ for pressure traces and $\Noisemeas=0 \; \si{\milli\liter}$, $\Noisemeas=1 \; \si{\milli\liter}$ and $\Noisemeas=10 \; \si{\milli\liter}$ for volume traces, respectively. We express the measurement errors by computing the signal-to-noise ratio as SNR$=\tfrac{\Noisemeas}{\mu}$, being $\mu = 100 \; \si{\mmHg}$ (equivalently, $\mu = 100 \; \si{\milli\liter}$) a reference value for the pressure (equivalently, volume) signal over time.

In the next sections, we compute the discrete $L^2$ relative error over the estimated parameters as follows:
\begin{equation*} 
E_{\text{L}^2}(\param) = \sqrt{ \dfrac{1}{\NumParams} \sum_{i=1}^{\NumParams} \left( \dfrac{\theta_{\text{i,MAP}} - \theta_{\text{i,exact}}}{\theta_{\text{i,exact}}} \right)^2}.
\end{equation*}

We also remark that we only use QoIs that can be measured non-invasively in patient-specific cases by means of standard techniques.
For instance, the systemic arterial pressure can be reconstructed from the pointwise values provided by a sphygmomanometer, while the time traces of atrial and ventricular volumes may be extracted from Cine Magnetic Resonance Imaging.
This allows for the clinical exploitation of the proposed approach.

We always perform 5 heartbeats to reach the limit cycle and we evaluate the time traces in the cost functional on the last heartbeat only, as reported in Section~\ref{sec:methods:OCP}.
We employ a sampling rate of $10^{-2}$ $\si{\second}$, that is comparable to the imaging time resolution \cite{Kellman2009}.
The target is always provided by the pressures and volumes computed by the $\modEMfom-\modCirc$ model. 

\subsection{Global sensitivity analysis}
\label{sec:GSA}

\begin{table}
    \begin{center}
        \begin{tabular}{ rrrl }
            \toprule
            Parameter & Baseline & Unit & Description \\
            \midrule
            $\aXB$         & 250 & $\si{\mega\pascal}$ & Cardiomyocytes contractility \\
            \bottomrule
        \end{tabular}
        \caption{Most relevant parameter of the $\modEMredANN$ model after sensitivity analysis and associated baseline value for the numerical tests.}
        \label{tab:paramsM}
    \end{center}
\end{table}

\begin{table}
    \begin{center}
        \begin{tabular}{ llrl }
            \toprule
            Parameter & Baseline & Unit & Description \\
            \midrule
            $\Vheart$                 & 417              & \si{\milli\liter}            & Initial blood volume in the cardiovascular system  \\
            $\EaRV$                   & 0.55             & \si{\mmHg \per \milli\liter} & RV active elastance \\
            $\EpLA$, $\EpRA$          & 0.15, 0.05       & \si{\mmHg \per \milli\liter} & LA/RA passive elastance \\
            $\RarSYS$,  $\RvnSYS$ & 0.64, 0.32 & \si{\mmHg \second \per \milli\liter}  & Systemic arterial/venous resistance \\
            \bottomrule
        \end{tabular}
        \caption{Most relevant parameters of the $\modCirc$ model after sensitivity analysis and associated baseline values for the numerical tests. LA: left atrium, RA: right atrium, RV: right ventricle.}
        \label{tab:paramsC}
    \end{center}
\end{table}

\begin{figure}
    \centering
    \includegraphics[width=1.0\textwidth]{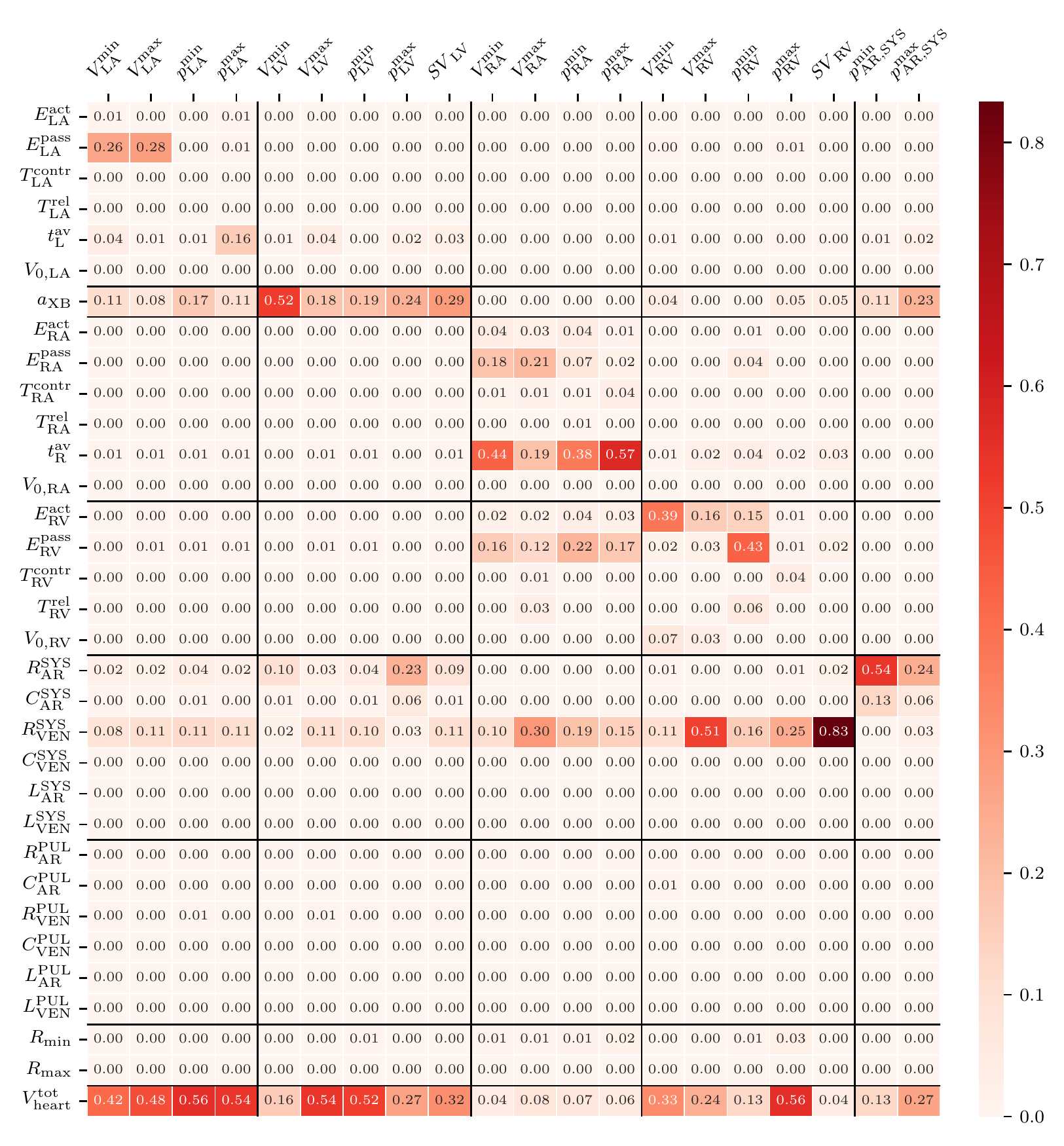}
    \caption{Total-effect Sobol indices computed by exploiting the $\modEMredANN-\modCirc$ model. For a detailed definition of the all model parameters and QoIs, we refer to Appendix~\ref{app:circulation}.}
    \label{fig:SA}
\end{figure}

As a preliminary step, we perform a global sensitivity analysis by employing the Saltelli's method \cite{saltelli2002making,homma1996importance}.
We compute both first-order and total-effect Sobol indices \cite{sobol1990sensitivity}.
The former quantifies how much varying a single parameter affects a specific QoI, while the latter also considers higher-order interactions among the parameters of the model.
For the sake of simplicity, we only depict total-effect Sobol indices in Figure~\ref{fig:SA}.
We report in Tables~\ref{tab:paramsM} and~\ref{tab:paramsC} the relevant parameters of the $\modEMredANN-\modCirc$ model, along with dimensional units, description and baseline values for the $\TLV$, $\Tall$ and $\TLVperturbed$ numerical tests.
As observed in our previous work \cite{Regazzoni2022EMROM}, heartbeat period, initial conditions in terms of blood pool volume in the heart, elastances of the cardiac chambers, timings of contraction and relaxation, atrioventricular delay and a couple of electrical components of the cardiocirculatory network play an effective role in determining pressures and volumes.
However, if volume traces over time are available, as in the test cases proposed in this paper, the heartbeat period along with the timings of contraction and relaxation of the different cardiac chambers can be directly observed and do not need to be estimated.
Regarding the $\modEMredANN$ model, contractility is certainly the most relevant parameter in determining the pressure-volume relationship.
All the relevant information on parameter bounds used to determine the MAP estimation and the values of the fixed parameters of the $\modEMredANN-\modCirc$ model are provided in Appendix~\ref{app:params}.

\subsection{$\TLV$ test case}
\label{sec:results:TLV}

We estimate $\param = [\aXB, \RarSYS, \Vheart]^T$ ($\NumParams$=3) with UQ starting from the time evolution of $\ParSYS$ and $\VLV$ by employing $\modEMredANN-\modCirc$ model.
We assume all the other parameters of the mathematical model to be known.
For the optimization process, we initialize $\aXB$, $\RarSYS$ and $\Vheart$ randomly in a uniform interval ranging from 50\% to 150\% of their exact values, in order to mimic a scenario in which there is no prior knowledge about $\param$.
We repeat the MAP estimation 10 times to leverage this random initialization.
We also test the robustness of our approach against different levels of SNR on the two time-dependent QoIs.

\begin{figure}
    \centering
    \includegraphics[width=0.7\textwidth]{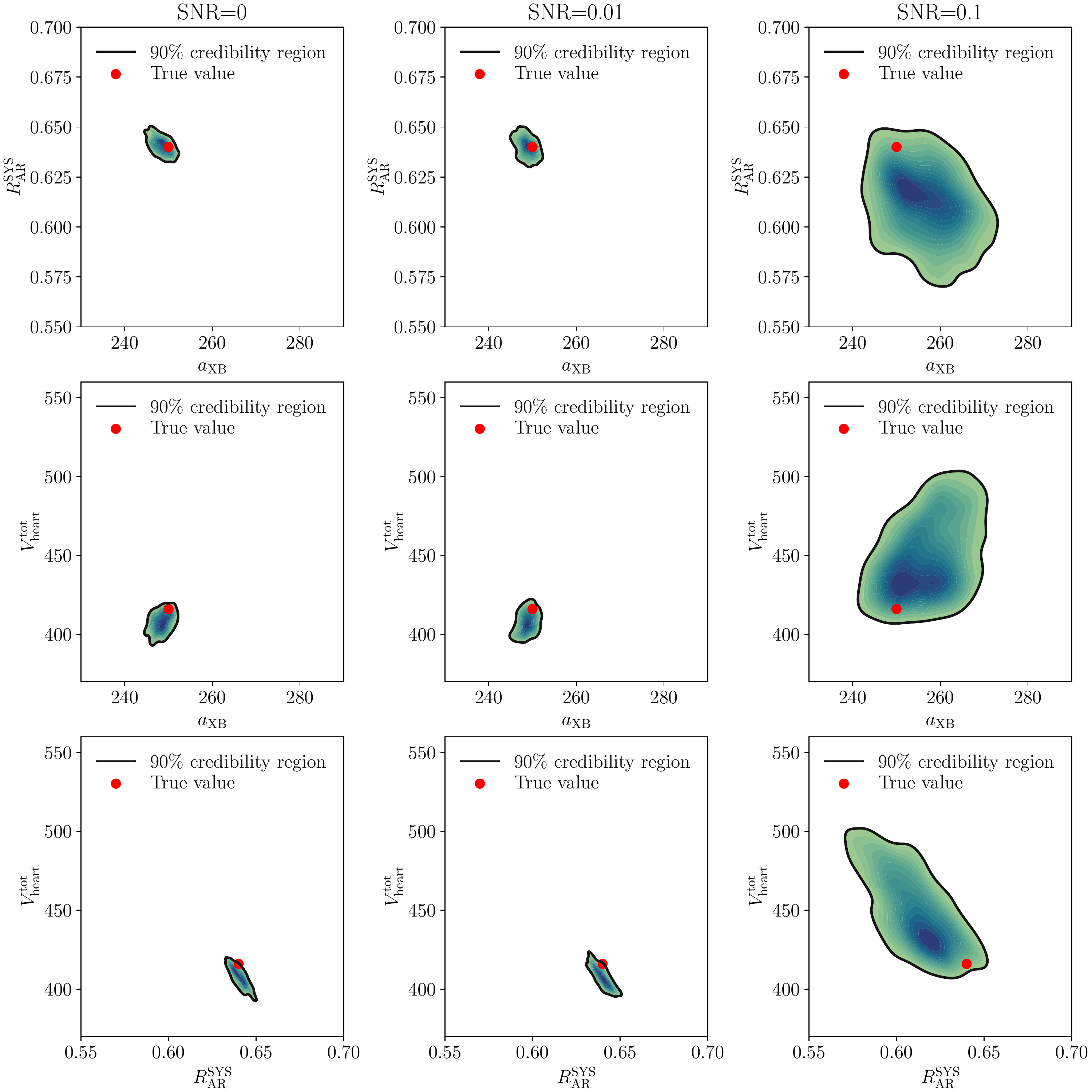}
    \caption{$\TLV$: posterior distribution $\mathbb{P}(\param | \mathbf{x})$ estimated by means of HMC. The red dots representing the exact values of $\aXB$, $\RarSYS$ and $\Vheart$ are included in the 90\% credibility regions.}
    \label{fig:TLVposterior}
\end{figure}

\begin{figure}
    \centering
    \includegraphics[width=1.0\textwidth]{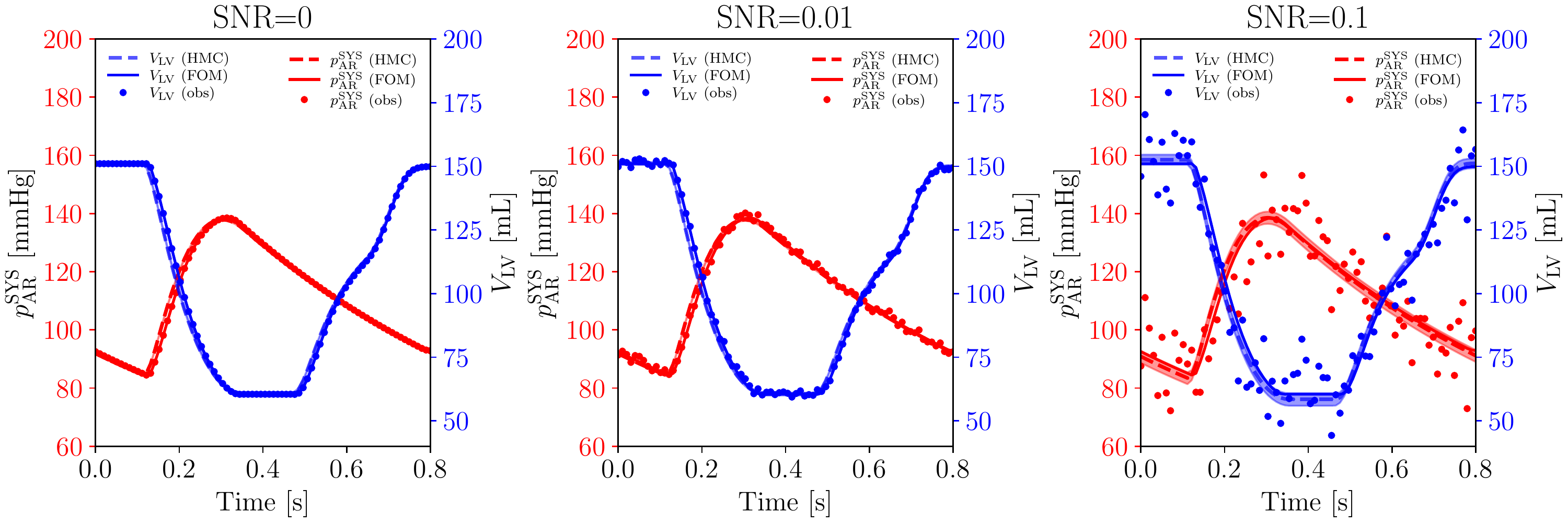}
    \caption{$\TLV$: $\ParSYS$ and $\VLV$ time traces for different levels of SNR. We show the FOM time traces, the corresponding noisy observations and the averaged HMC time traces endowed with the $5^{\text{th}}$ and $95^{\text{th}}$ percentiles.}
    \label{fig:TLVpressurevolume}
\end{figure}

We report in Table~\ref{tab:TLVMAP} some statistics regarding these 10 runs, such as the value of $E_{\text{L}^2}(\param)$ computed by averaging the 10 estimated values of the parameters, the number of L-BFGS iterations, the final loss function value in the best and worst case scenarios, and the total computational time.
Then, we perform HMC with the same three values of SNR by initializing the algorithm with the value of $\param$ averaged over the 10 MAP estimations. We report in Table~\ref{tab:TLVHMC} the total execution time while increasing the noise level.
We show in Table~\ref{tab:TLVHMCtruemeanstd} the computed mean and standard deviation of $\param$ by considering the HMC trajectory after the warmup period for different values of SNR.
We display in Figure~\ref{fig:TLVposterior} the posterior distribution $\mathbb{P}(\param | \mathbf{x})$ for the SNR going from 0 to 0.1, which always contains the true value of $\param$ in the $90\%$ credibility region.
We show in Figure~\ref{fig:TLVpressurevolume} the QoIs, i.e. the time evolution of $\ParSYS$ and $\VLV$.
The FOM solution is mostly contained inside the $5^{\text{th}}$ and $95^{\text{th}}$ percentiles of the HMC estimations from low-to-high noise levels.

\begin{table}[t!]
    \begin{center}
        \hspace*{-0.25cm}
        \begin{tabular}{ l c cc cc c}
            \toprule
            SNR & $E_{\text{L}^2}(\param)$ & \multicolumn{2}{c}{L-BFGS iterations} & \multicolumn{2}{c}{Loss function $J(\param)$} & Total execution time [s] \\
            \midrule
             & Average & Best & Worst & Best & Worst & \\
            0 & 7.23e-03 & 20 & 48 & 5.32e-05 & 2.93e-04 & 28 minutes \\
            0.01 & 3.87e-02 & 17 & 50 & 2.30e-04 & 1.19e-03 & 26 minutes \\
            0.1 & 7.46e-02 & 19 & 25 & 1.64e-02 & 1.70e-02 & 21 minutes \\
            \bottomrule
        \end{tabular}
        \caption{$\TLV$: MAP estimation with different values of SNR. The $E_{\text{L}^2}(\param)$ is obtained by averaging 10 MAP estimations with a random initialization of $\param$ in a uniform interval. A single run approximately takes 2-3 minutes of computations on average.}
        \label{tab:TLVMAP}
    \end{center}
\end{table}

\begin{table}[t!]
    \begin{center}
        \begin{tabular}{ ll }
            \toprule
            SNR & Total execution time [s] \\
            \midrule
            0 & 2 hours and 29 minutes \\
            0.01 & 2 hours and 6 minutes \\
            0.1 & 1 hours and 57 minutes \\
            \bottomrule
        \end{tabular}
        \caption{$\TLV$: execution times for HMC with different values of SNR.}
        \label{tab:TLVHMC}
    \end{center}
\end{table}

\begin{table}[t!]
    \begin{center}
        \begin{tabular}{ l|rrrr }
            \toprule
            Parameter & Reference & SNR = 0 & SNR = 0.1 & SNR = 0.01 \\
            \midrule
            $\aXB$    & 250.00 & 248.46 $\pm$ 3.34 & 248.78 $\pm$ 3.07 & 256.66 $\pm$ 13.64 \\
            $\RarSYS$ & 0.64 & 0.64 $\pm$ 0.01& 0.64 $\pm$ 0.01& 0.61 $\pm$ 0.03\\ 
            $\Vheart$ & 417 & 407 $\pm$ 12& 408 $\pm$ 12 & 448 $\pm$ 43\\
            \bottomrule
        \end{tabular}
        \caption{$\TLV$: mean plus/minus two standard deviations associated to the estimated parameters during HMC for different values of SNR. We refer to Tables~\ref{tab:paramsM} and~\ref{tab:paramsC} for the ground truth.}
        \label{tab:TLVHMCtruemeanstd}
    \end{center}
\end{table}

\subsection{$\Tall$ test case}
\label{sec:results:Tall}

We consider a complicated framework in which, given non-invasive time traces of $\ParSYS$, $\VLA$, $\VLV$, $\VRA$ and $\VRV$, we estimate $\param = [\aXB, \EaRV, \EpLA, \EpRA, \RarSYS, \RvnSYS, \Vheart]^T$ ($\NumParams$=7).
Again, we assume all the other parameters of the mathematical model to be known.
To overcome the random initialization of $\param$ in the 50\% to 150\% range discussed in Section~\ref{sec:results:TLV}, we still run the MAP estimation 10 times. Indeed, given the robustness of the L-BFGS method with respect to its initialization, the number of runs does not have to vary with the number of parameters and QoIs.

\begin{figure}
    \centering
    \includegraphics[width=0.6\textwidth]{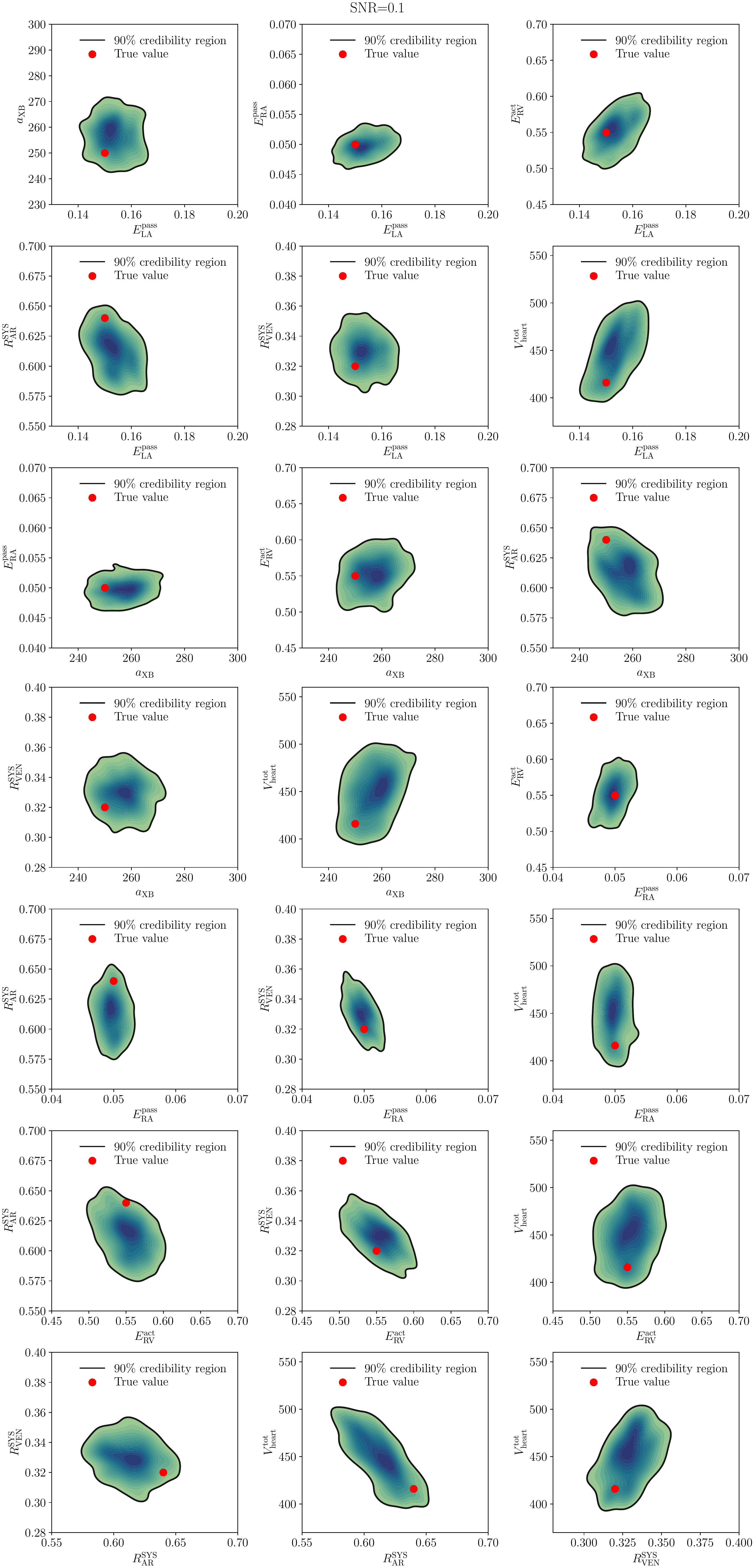}
    \caption{$\Tall$: posterior distribution $\mathbb{P}(\param | \mathbf{x})$ estimated by means of HMC for SNR=0.1. The red dots representing the exact values of $\aXB$, $\EaRV$, $\EpLA$, $\EpRA$, $\RarSYS$, $\RvnSYS$ and $\Vheart$ are contained in the 90\% credibility regions.}
    \label{fig:Tallposterior}
\end{figure}

The statistics regarding the MAP estimations are reported in Table~\ref{tab:TallMAP}.
We collect the execution times related to the HMC runs in Table~\ref{tab:TallHMC}, while we depict the posterior distribution for the case SNR=0.1 in Figure~\ref{fig:Tallposterior}.
In Table~\ref{tab:TallHMCtruemeanstd} we show the computed mean and standard deviation of all the parameters $\param$ by considering the samples coming from HMC for different noise levels.

\subsection{$\TLVperturbed$ test case}
\label{sec:results:TLVperturbed}

As in $\TLV$, we perform the MAP estimation for $\param = [\aXB, \RarSYS, \Vheart]^T$ ($\NumParams$=3) starting from the time evolution of $\ParSYS$ and $\VLV$ by employing the $\modEMredANN-\modCirc$ model.
However, here, we assume that other relevant parameters of the mathematical model ($\EaRV, \EpLA, \EpRA$ and $\RvnSYS$) are not set to the ground truth, as it is supposed to be in realistic scenarios.
For the sake of simplicity, these four parameters may be uniformly assumed to be around $5\%$, $10\%$ or $20\%$ of their exact values in the 10 MAP estimation runs.
We report the numerical results in Table~\ref{tab:TLVperturbedMAP} for two different values of SNR, 0 and 0.05 respectively.

\begin{table}[t!]
    \begin{center}
        \hspace*{-0.25cm}
        \begin{tabular}{ c c cc cc c}
            \toprule
            SNR & $E_{\text{L}^2}(\param)$ & \multicolumn{2}{c}{L-BFGS iterations} & \multicolumn{2}{c}{Loss function $J(\param)$} & Total execution time [s] \\
            \midrule
             & Average & Best & Worst & Best & Worst & \\
            0 & 3.41e-02 & 19 & 50 & 7.77e-05 & 4.83e-03 & 28 minutes \\
            0.01 & 3.55e-02 & 42 & 75 & 6.05e-04 & 1.86e-03 & 34 minutes \\
            0.1 & 4.15e-02 & 19 & 50 & 5.17e-02 & 5.49e-02 & 28 minutes \\
            \bottomrule
        \end{tabular}
        \caption{$\Tall$: MAP estimation with different values of SNR. The $E_{\text{L}^2}(\param)$ is obtained by averaging 10 MAP estimations with a random initialization of $\param$ in a uniform interval. A single run approximately takes 2-4 minutes of computations on average.}
        \label{tab:TallMAP}
    \end{center}
\end{table}

\begin{table}[t!]
    \begin{center}
        \begin{tabular}{ ll }
            \toprule
            SNR & Total execution time [s] \\
            \midrule
            0 & 4 hours and 50 minutes \\
            0.01 & 4 hours and 35 minutes \\
            0.1 & 4 hours and 13 minutes \\
            \bottomrule
        \end{tabular}
        \caption{$\Tall$: execution times for HMC with different values of SNR.}
        \label{tab:TallHMC}
    \end{center}
\end{table}

\begin{table}[t!]
    \begin{center}
        \hspace*{-0.5cm}
        \begin{tabular}{ l|rrrr }
            \toprule
            Parameter & Reference & SNR = 0 & SNR = 0.1 & SNR = 0.01 \\
            \midrule
            $\aXB$  & 250.00 & 247.89 $\pm$ 2.74 & 248.91 $\pm$ 3.17 & 256.93 $\pm$ 12.64 \\
            $\EaRV$ & 0.550 & 0.548 $\pm$ 0.003 & 0.545 $\pm$ 0.007 & 0.553 $\pm$ 0.044 \\ 
            $\EpLA$ & 0.150 & 0.148 $\pm$ 0.002 & 0.147 $\pm$ 0.003 & 0.154 $\pm$ 0.011 \\
            $\EpRA$ & 0.050 & 0.0497 $\pm$ 0.0001 & 0.0499 $\pm$ 0.0004 & 0.0498 $\pm$ 0.0032 \\
            $\RarSYS$ & 0.640 & 0.636 $\pm$ 0.006 & 0.637 $\pm$ 0.008 & 0.613 $\pm$ 0.033 \\
            $\RvnSYS$ & 0.320 & 0.323 $\pm$ 0.001 & 0.321 $\pm$ 0.003 & 0.329 $\pm$ 0.022 \\
            $\Vheart$ & 417 & 422 $\pm$ 7 & 415 $\pm$ 12 & 449 $\pm$ 48 \\
            \bottomrule
        \end{tabular}
        \caption{$\Tall$: mean plus/minus two standard deviations associated to the estimated parameters during HMC for different values of SNR. We refer to Tables~\ref{tab:paramsM} and~\ref{tab:paramsC} for the ground truth.}
        \label{tab:TallHMCtruemeanstd}
    \end{center}
\end{table}

\begin{table}[t!]
    \begin{center}
        \hspace*{-0.75cm}
        \begin{tabular}{ c c c cc cc c}
            \toprule
            SNR & Perturbation & $E_{\text{L}^2}(\param)$ & \multicolumn{2}{c}{L-BFGS iterations} & \multicolumn{2}{c}{Loss function $J(\param)$} & Total execution time [s] \\
            \midrule
             & & Average & Best & Worst & Best & Worst & \\
            0 & 5 $\%$ & 6.71e-03 & 19 & 41 & 5.37e-05 & 1.19e-04 & 24 minutes \\
            0 & 10 $\%$ & 3.47e-03 & 16 & 50 & 5.39e-05 & 2.01e-04 & 24 minutes \\
            0 & 20 $\%$ & 2.70e-03 & 14 & 45 & 5.44e-05 & 4.08e-04 & 24 minutes \\
            \midrule
            0.05 & 5 $\%$ & 1.86e-02 & 21 & 47 & 4.19e-03 & 4.52e-03 & 25 minutes \\
            0.05 & 10 $\%$ & 4.35e-02 & 17 & 57 & 4.19e-03 & 4.74e-03 & 26 minutes \\
            0.05 & 20 $\%$ & 2.57e-02 & 20 & 47 & 4.19e-03 & 4.79e-03 & 22 minutes \\            
            \bottomrule
        \end{tabular}
        \caption{$\TLVperturbed$: MAP estimation with different values of SNR and perturbation of significant parameters of the $\modEMredANN-\modCirc$ model. The $E_{\text{L}^2}(\param)$ is obtained by averaging 10 MAP estimations with a random initialization of $\param$ in a uniform interval. A single run approximately takes 2-3 minutes of computations on average.}
        \label{tab:TLVperturbedMAP}
    \end{center}
\end{table}

	\section{Discussion}
\label{sec:discussion}

From the results reported in Tables~\ref{tab:TLVMAP},~\ref{tab:TallMAP} and~\ref{tab:TLVperturbedMAP}, we conclude that the MAP estimation process is robust from low-to-high values of SNR.
Indeed, even though $E_{\text{L}^2}(\param)$ increases with SNR, its maximum value, obtained by averaging the final values of $\param$ over the 10 runs, is always under $8\%$, and in some cases lower than $1\%$.
In particular, for the $\Tall$ test case (see Table~\ref{tab:TallMAP}), the values of $E_{\text{L}^2}(\param)$ for SNR=0 and SNR=0.01 are very similar.
This happens when the ROM approximation error is comparable to the measurement error.
Indeed, a cancellation phenomenon between the two sources of error may occur, leading to similar performances with respect to the case where there is no measurement error (SNR=0).
Furthermore, the final value of the loss function decreases as the noise increases, because it is intrinsically more challenging to fit the QoIs for higher values of SNR.
From Table~\ref{tab:TLVperturbedMAP}, we see that $E_{\text{L}^2}(\param)$ remains small and comparable when we consider increasing perturbations on influential parameters for the cardiac function that are not estimated during the procedure.
This is motivated by the fact that the two observed QoIs, namely $\ParSYS$ and $\VLV$, mostly depend on $\aXB$, $\RarSYS$, $\Vheart$.
This can be observed from the global sensitivity analysis reported in Figure~\ref{fig:SA}.

Finally, thanks to the use of matrix-free methods, automatic differentiation and automatic vectorization, the 10 runs for MAP estimation always require a small total computational time, which is approximately equal to 30 minutes.
This is a remarkable result, especially given the use of just one CPU in a standard laptop.
Furthermore, the computational cost is comparable among the different levels of noise and the 10 runs are embarrassingly parallelizable.

The HMC method initialized with the MAP estimation proves to be very efficient and robust as well.
By looking at Tables~\ref{tab:TLVHMC} and~\ref{tab:TallHMC}, we see that the total execution time is always under 5 hours and decreases as the SNR increases.
This is justified by the fact that the trajectory associated to $\param$ may evolve more freely (consequently faster) in a high noise level scenario.
From Figures~\ref{fig:TLVposterior} and~\ref{fig:Tallposterior}, we notice that, as expected, the $90\%$ credibility regions of the estimated posterior distribution $\mathbb{P}(\param | \mathbf{x})$ get bigger as the SNR increases.
Moreover, they always contain the exact values for $\param$.
This is also quantitatively confirmed by the statistics reported in Tables~\ref{tab:TLVHMCtruemeanstd} and~\ref{tab:TallHMCtruemeanstd}, for all the estimated parameters and noise levels.
This Bayesian framework also allows to capture correlations among the estimated parameters.
For instance, for the range of values that have been considered in this paper, $\RarSYS$ is negatively correlated with $\aXB$ and $\Vheart$, $\RvnSYS$ and $\EpLA$ are positively correlated with $\Vheart$, while varying $\Vheart$ or $\aXB$ does not affect the value of $\EpRA$.
	\section{Conclusions}
\label{sec:conclusions}

We employ Bayesian statistics to perform parameter estimation with UQ for the cardiac function.
We effectively combine MAP estimation and HMC in different in silico test cases and we use an ANN-based surrogate model of 3D cardiac electromechanics coupled with a 0D model of the cardiocirculatory network, thus allowing for a massive reduction of the computational cost associated with the parameter estimation procedure.
We show that we can robustly estimate several parameters, associated to both cardiac mechanics and cardiovascular hemodynamics models, starting from noisy non-invasive pressure-volume measurements and a random initialization of the model parameters.
We always reach relative $L^2$ errors under $8\%$ on the estimated parameters and meaningful posterior distributions in less than 5 hours of computations by using a single CPU of a standard laptop.
We incorporate the presence of both model approximation and measurement errors in our method.
We always consider a number of estimated parameters that is less than the number of observed time traces in the cardiovascular system.
This is certainly important in realistic scenarios, where clinical data are generally both noisy and lacking.

As future developments, we aim at extending the ANN-based surrogate model to account for different cardiac chambers and multiple 3D model parameters, as the LV, along with contractility, are uniquely considered in the present work.
Furthermore, we also want to consider clinical data, rather than in silico numerical simulations coming from a 3D-0D model.
Finally, we remark that the software libraries on which this method is based are natively compatible with Graphical Processing Units and Tensor Processing Units. Parallel computing is also supported for all the hardware architectures.
This may significantly reduce the total computational costs associated to MAP estimation and HMC.

	\section*{Acknowledgements}
This project has been funded by the Italian Ministry of University and Research (MIUR) within the PRIN (Research projects of relevant national interest 2017 ``Modeling the heart across the scales: from cardiac cells to the whole organ'' Grant Registration number 2017AXL54F).

\begin{center}
	\raisebox{-.5\height}{\includegraphics[width=.15\textwidth]{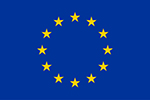}}
	\hspace{2cm}
	\raisebox{-.5\height}{\includegraphics[width=.15\textwidth]{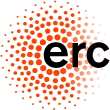}}
\end{center}

	\begin{appendices}
		\section{Circulation model}
\label{app:circulation}


We model fluid dynamics within the entire cardiovascular system by means of a 0D closed-loop circulation model $\modCirc$ proposed in \cite{Regazzoni2022} and inspired by \cite{Hirschvogel2017}.
We employ time-varying elastances for the four cardiac chambers, non-ideal diodes for the four cardiac valves and resistance-inductance-capacitance (RLC) circuits for the systemic and pulmonary circulations, which in turn can be divided into the arterial and venous branches.
The 0D closed-loop circulation model reads:
\begin{subequations}
\begin{empheq}[left={\empheqlbrace\,}]{align}
\dfrac{d \VLA(t)}{d t} & = \QvnPUL(t) - \QMV(t), \label{eqn: circulation_part1: VLA} \\
\dfrac{d \VLV(t)}{d t} & = \QMV(t)    - \QAV(t), \label{eqn: circulation_part1: VLV} \\
\dfrac{d \VRA(t)}{d t} & = \QvnSYS(t) - \QTV(t), \label{eqn: circulation_part1: VRA} \\
\dfrac{d \VRV(t)}{d t} & = \QTV(t)    - \QPV(t), \label{eqn: circulation_part1: VRV} \\
\CarSYS \dfrac{d \ParSYS(t)}{d t} & = \QAV(t)    - \QarSYS(t), \label{eqn: circulation_part1: ParSYS} \\
\CvnSYS \dfrac{d \PvnSYS(t)}{d t} & = \QarSYS(t) - \QvnSYS(t), \label{eqn: circulation_part1: PvnSYS} \\
\CarPUL \dfrac{d \ParPUL(t)}{d t} & = \QPV(t)    - \QarPUL(t), \label{eqn: circulation_part1: ParPUL} \\
\CvnPUL \dfrac{d \PvnPUL(t)}{d t} & = \QarPUL(t) - \QvnPUL(t), \label{eqn: circulation_part1: PvnPUL} \\
\dfrac{\LarSYS}{\RarSYS} \dfrac{d \QarSYS(t)}{d t} & = - \QarSYS(t) - \dfrac{\PvnSYS(t) - \ParSYS(t)}{\RarSYS}, \label{eqn: circulation_part1: QarSYS} \\
\dfrac{\LvnSYS}{\RvnSYS} \dfrac{d \QvnSYS(t)}{d t} & = - \QvnSYS(t) - \dfrac{\PRA(t)    - \PvnSYS(t)}{\RvnSYS}, \label{eqn: circulation_part1: QvnSYS} \\
\dfrac{\LarPUL}{\RarPUL} \dfrac{d \QarPUL(t)}{d t} & = - \QarPUL(t) - \dfrac{\PvnPUL(t) - \ParPUL(t)}{\RarPUL}, \label{eqn: circulation_part1: QarPUL} \\
\dfrac{\LvnPUL}{\RvnPUL} \dfrac{d \QvnPUL(t)}{d t} & = - \QvnPUL(t) - \dfrac{\PLA(t)    - \PvnPUL(t)}{\RvnPUL}, \label{eqn: circulation_part1: QvnPUL}
\end{empheq}
\label{eqn: circulation_part1}
\end{subequations}
with $t \in (0, T]$, where:
\begin{subequations}
\begin{align}
\PLV(t) & = \ELV(t) \left(\VLV(t) - \VnLV \right), \label{eqn: circulation_part2: PLV} \\
\PLA(t) & = \ELA(t) \left(\VLA(t) - \VnLA \right), \label{eqn: circulation_part2: PLA} \\
\PRV(t) & = \ERV(t) \left(\VRV(t) - \VnRV \right), \label{eqn: circulation_part2: PRV} \\
\PRA(t) & = \ERA(t) \left(\VRA(t) - \VnRA \right), \label{eqn: circulation_part2: PRA} \\
\QMV(t) & = \dfrac{\PLA(t) -\PLV(t)    }{\RMV(\PLA(t), \PLV(t))}, \label{eqn: circulation_part2: QMV} \\
\QAV(t) & = \dfrac{\PLV(t) -\ParSYS(t) }{\RAV(\PLV(t), \ParSYS(t))}, \label{eqn: circulation_part2: QAV} \\
\QTV(t) & = \dfrac{\PRA(t) -\PRV(t)    }{\RTV(\PRA(t), \PRV(t))}, \label{eqn: circulation_part2: QTV} \\
\QPV(t) & = \dfrac{\PRV(t) -\ParPUL(t) }{\RPV(\PRV(t), \ParPUL(t))}, \label{eqn: circulation_part2: QPV}
\end{align}
\label{eqn: circulation_part2}
\end{subequations}
\noindent with $ t \in (0, T]$.
$ \PLA(t) $, $ \PRA(t) $, $ \PLV(t) $, $ \PRV(t) $, $ \VLA(t) $, $ \VRA(t) $, $ \VLV(t) $ and $ \VRV(t) $ refer to pressures and volumes in the four cardiac chambers.
We remark that, if we consider the coupling with $\modEMfom$, $\modEMredANN$ or $\modEMredemulator$ electromechanical models for the LV, then Equation~\eqref{eqn: circulation_part2: PLV} is not solved \cite{Regazzoni2022}.
The same approach can be seamlessly extended to the atria and right ventricle as well \cite{Fedele2022,Piersanti2022}.
The fluxes accross the cardiac valves are expressed by $ \QMV(t) $, $ \QAV(t) $, $ \QTV(t) $ and $ \QPV(t) $.
Arterial and venous systemic circulation is described by means of $ \ParSYS(t) $, $ \QarSYS(t) $, $ \PvnSYS(t) $ and $ \QvnSYS(t) $, while $ \ParPUL(t) $, $ \QarPUL(t) $, $ \PvnPUL(t) $ and $ \QvnPUL(t) $ define the behavior of the arterial and venous pulmonary circulation.
Time-varying elastances $ \ELA(t) $, $ \ELV(t) $, $ \ERA(t) $, $ \ERV(t) $ are analytically prescribed functions that define the pressure-volume relationship of atria and ventricles. They range from $ \EpLA $, $ \EpLV $, $ \EpRA $, $ \EpRV $ -- when the chambers are at rest -- to $ (\EpLA+\EaMaxLA) $, $ (\EpLV+\EaMaxLV) $, $ (\EpRA+\EaMaxRA) $, $ (\EpRV+\EaMaxRV) $ -- when the chambers are fully contracted. Specifically:
\begin{equation*}
	\actshape(t;\tC,\tR,\TC,\TR) = \begin{cases}
		\frac{1}{2} \left[1 - \cos\left( \frac{\pi}{\TC}\operatorname{mod}(t - \tC, \THB) \right) \right]
		&
		\text{if $0 \leq \operatorname{mod}(t - \tC, \THB) < \TC $},
		\\[10pt]
		\frac{1}{2} \left[1 + \cos\left( \frac{\pi}{\TR}\operatorname{mod}(t - \tR, \THB) \right) \right]
		&
		\text{if $0 \leq \operatorname{mod}(t - \tR, \THB) < \TR $},
		\\[10pt]
		0
		&
		\text{otherwise,}
	\end{cases}
\end{equation*}
where $\tC$ and $\tR$ are the contraction and relaxation initial times, respectively and $\TC$ and $\TR$ are the contraction and relaxation durations, respectively. We remark that $\tR = \tC + \TC$ for consistency.
Then, we write the different time-varying elastances as follows:
\begin{equation*}
	\Estar(t) = \Epstar + \EaMaxstar \actshape(t, \tCstar,\tRstar,\TCstar,\TRstar),
	\qquad
	\text{for} \, i \in \{\text{LA}, \text{LV}, \text{RA}, \text{RV}\}.
\end{equation*}
Again, we underline that when $\modEMfom$, $\modEMredANN$ or $\modEMredemulator$ electromechanical models are coupled to the circulation model $\modCirc$, the time-varying elastances are replaced by a mechanical activation model \cite{regazzoni2020biophysically,Salvador2020}.

Finally, we use diodes to model the four cardiac valves:
\begin{equation*}
	R_{\text{i}}(p_1, p_2) =
	\begin{cases}
		\Rmin, & p_1 < p_2 \\
		\Rmax, & p_1 \geq p_2 \\
	\end{cases}
	\;\;\; \text{for} \;\;\; i \in \{\text{MV}, \text{AV}, \text{TV}, \text{PV}\},
\end{equation*}
where $p_1$ and $p_2$ denote the pressures ahead and behind the valve leaflets with respect to the flow direction, whereas $\Rmin$ and $\Rmax$ are the minimum and maximum resistance of the valves.

	
		\section{Parameter values}
\label{app:params}

We vary the parameters of Tables~\ref{tab:paramsM} and~\ref{tab:paramsC} in suitable ranges during the bounded constrained optimization process.
In particular, we retrieve the values for relevant electrical components from \cite{Hirschvogel2017, Klingensmith2008, Saouti2010}, that is $\RarSYS \in [0.54, 1.2]$ \si{\mmHg \second \per \milli\liter} and $\RvnSYS \in [0.18, 0.4]$ \si{\mmHg \second \per \milli\liter}.
We vary all the properties related to time-varying elastances from 50\% to 150\% of their baseline values \cite{Regazzoni2022}.
We consider $\aXB \in [80, 320]$ $\si{\mega\pascal}$ for $\modEMredANN$ model \cite{Regazzoni2022EMROM}. 
Finally, we define $\Vheart \in [200, 600]$ \si{\milli\liter} \cite{Lemmens2006, Nadler1962}.
In this way, we can also simulate pathological scenarios, such as mild dilated cardiomyopathy with the possible presence of ischemia \cite{Salvador2021,Salvador2021MEF}.
The other parameters of the $\modCirc$ model are fixed to the values reported in Table~\ref{tab:paramcirc}.

\begin{table}[h!]
	\centering
	\begin{tabular}{lrr|lrr}
		\toprule
		Variable & Value & Unit & Variable & Value & Unit \\
		\midrule
		\multicolumn{3}{l|}{\textbf{External circulation}} & \multicolumn{3}{l}{\textbf{Cardiac chambers}} \\
		$\RarPUL$    & 0.032      & \si{\mmHg \second \per \milli\liter}   &
		$\EpRV$      & 0.05       & \si{\mmHg \per \milli\liter}   \\
		$\RvnPUL$    & 0.035      & \si{\mmHg \second \per \milli\liter}   &
		$\EaMaxLA$   & 0.07       & \si{\mmHg \per \milli\liter}   \\
		$\CarSYS$    & 1.2          & \si{\milli\liter \per \mmHg}   &
		$\EaMaxRA$   & 0.20       & \si{\mmHg \per \milli\liter}   \\
		$\CarPUL$    & 10         & \si{\milli\liter \per \mmHg}   &
		 &  & \\
		$\CvnSYS$    & 60         & \si{\milli\liter \per \mmHg}   &
		$\VnLA$      & 4          & \si{\milli\liter}   \\
		$\CvnPUL$    & 16         & \si{\milli\liter \per \mmHg}   &
		$\VnRA$      & 4          & \si{\milli\liter}   \\
		$\LarSYS$    & \num{5e-3}         & \si{\mmHg \second\squared \per \milli\liter}  &
		$\VnRV$      & 16         & \si{\milli\liter}   \\
		$\LarPUL$    & \num{5e-4}         & \si{\mmHg \second\squared \per \milli\liter}  &
		$\tC_\text{LA}$ & 0.64 & \si{\second} \\
		$\LvnSYS$ & \num{5e-4} & \si{\mmHg \second\squared \per \milli\liter}  &
		$\TC_\text{LA}$ & 0.12 & \si{\second} \\
		$\LvnPUL$ & \num{5e-4} & \si{\mmHg \second\squared \per \milli\liter}  &
		$\TR_\text{LA}$ & 0.64 & \si{\second} \\
		\multicolumn{3}{l|}{\textbf{Cardiac valves}} &
		$\tC_\text{RA}$ & 0.69 & \si{\second} \\
		$\Rmin$      & 0.0075       & \si{\mmHg \second \per \milli\liter} &
		$\TC_\text{RA}$ & 0.08 & \si{\second} \\
		$\Rmax$      & 75000      & \si{\mmHg \second \per \milli\liter} &
		$\TR_\text{RA}$ & 0.56 & \si{\second} \\
		& & &
		$\tC_\text{RV}$ & 0.04 & \si{\second} \\
		& & &
		$\TC_\text{RV}$ & 0.20 & \si{\second} \\
		& & &
		$\TR_\text{RV}$ & 0.32 & \si{\second} \\
		\bottomrule
	\end{tabular}
	\caption{Fixed parameters of the $\modCirc$ model. We consider a heartbeat period $\THB = \SI{0.8}{\second}$.}
	\label{tab:paramcirc}
\end{table}

	\end{appendices}
    \newpage
	\printbibliography

\end{document}